\documentclass[12pt]{amsart}
\usepackage{pifont}  
\usepackage{longtable}
\usepackage{amsfonts,amssymb,amscd,amsmath,epsfig}
\usepackage{latexsym}
\usepackage{mathrsfs}
\usepackage{tocvsec2}
\usepackage[]{youngtab}
\usepackage[%dvipdfm,
hyperindex=true,pagebackref=true,bookmarks=true,
colorlinks=true,linkcolor=blue,citecolor=red]{hyperref}
\let\oin\o
\let\cin\c
\let\vin\v
\let\lin\l

\begin{document}

%\oin, \cin{a}, \vin{a}, \lin

%MACOS FOR LECTURES ON DAHA
\renewcommand{\tilde}{\widetilde}
\renewcommand{\hat}{\widehat}

\newcommand{\BR}{{\mathbb R}}
\newcommand{\BQ}{{\mathbb Q}}
\newcommand{\BC}{{\mathbb C}}
\newcommand{\BP}{{\mathbb P}}
\newcommand{\BZ}{{\mathbb Z}}
\newcommand{\BN}{{\mathbb N}}
\newcommand{\BS}{{\mathbb S}}

\newcommand{\cH}{{\mathcal H}}
\newcommand{\cA}{{\mathcal A}}
\newcommand{\cB}{{\mathcal B}}
\newcommand{\ccF}{{\mathfrak F}}
\newcommand{\cD}{{\mathcal D}}
\newcommand{\cL}{{\mathcal L}}
\newcommand{\cF}{{\mathcal F}}
\newcommand{\cP}{{\mathcal P}}
\newcommand{\cX}{{\mathcal X}}
\newcommand{\cY}{{\mathcal Y}}
\newcommand{\cS}{{\mathcal S}}
\newcommand{\cSol}{\hbox{$\mathcal Sol$}}
\newcommand{\cT}{\hbox{$\mathcal T$}}

\newcommand{\Z}{{\mathbb Z}}
\newcommand{\Q}{{\mathbb Q}}
\newcommand{\N}{{\mathbb N}}
\newcommand{\C}{{\mathbb C}}
\newcommand{\R}{{\mathbb R}}
\newcommand{\X}{{\mathbb X}}
\newcommand{\Y}{{\mathbb Y}}

\newcommand{\CH}{{\mathcal H}}
\newcommand{\CA}{{\mathcal A}}

\def\HH{\mbox{${\mathcal H}$\kern-5.2pt${\mathcal H}$}}

\newcommand{\binomial}[2]{\genfrac{(}{)}{0pt}{}{ #1 }{ #2 }}
\newcommand{\qbinomial}[2]{\genfrac{[}{]}{0pt}{}{ #1 }{ #2 }_q }
\newcommand{\qbinom}[3]{\genfrac{[}{]}{0pt}{}{ #1 }{ #2 }_{ #3 } }

%%SPECIAL SEC 1.0

\def\der{\partial}
\def\tensor{\otimes}
\def\gam{\gamma} \def\Gam{\Gamma}
\def\del{\delta} \def\Del{\Delta}
\def\kap{\kappa}
\def\lam{\lambda} \def\Lam{\Lambda}
\def\Comp{{\mathbb C}}
\def\sM{{\mathcal M}}

\newtheorem{theorem}{Theorem}[section]
\newtheorem{maintheorem}[theorem]{Main Theorem}
\newtheorem{proposition}[theorem]{Proposition}
\newtheorem{definition}[theorem]{Definition}
\newtheorem{lemma}[theorem]{Lemma}
\newtheorem{corollary}[theorem]{Corollary}
\newtheorem{notation}[theorem]{Notation}
\newtheorem{remark}[theorem]{Remark}
\newtheorem{example}[theorem]{Example}

\newtheorem{theorem }{Theorem}[section]
\newtheorem{maintheorem }[theorem]{Main Theorem}
\newtheorem{proposition }[theorem]{Proposition}
\newtheorem{definition }[theorem]{Definition}
\newtheorem{lemma }[theorem]{Lemma}
\newtheorem{corollary }[theorem]{Corollary}
\newtheorem{notation }[theorem]{Notation}
\newtheorem{remark }[theorem]{Remark}
\newtheorem{example }[theorem]{Example}

\newtheorem{ maintheorem }[theorem]{Main Theorem}
\newtheorem{ theorem}{Theorem}[section]
\newtheorem{ proposition}[theorem]{Proposition}
\newtheorem{ definition}[theorem]{Definition}
\newtheorem{ lemma}[theorem]{Lemma}
\newtheorem{ corollary}[theorem]{Corollary}
\newtheorem{ notation}[theorem]{Notation}
\newtheorem{ remark}[theorem]{Remark}
\newtheorem{ example}[theorem]{Example}

\newtheorem{thm}{Theorem}[section]
\newtheorem{prop}[thm]{Proposition}
\newtheorem{lem}[thm]{Lemma}
\newtheorem{cor}[thm]{Corollary}
\newtheorem{conj}[thm]{Conjecture}
\newtheorem{con}[thm]{Conjecture}
\newtheorem{dfn}[thm]{Definition}
\newtheorem{df}[thm]{Definition}
 \newcommand{\rem}{{\bf Comment.\ }}
 \newcommand{\rmk}{{\bf Comment.\ }}
 \newcommand{\exmp}{{\bf Example.\ }}
 \newcommand{\ex}{{\bf Example.\ }}
 \newcommand{\prob}{{\bf Problem.\ }}

\newtheorem{note}{Note} 
\renewcommand{\thenote}{}
\newtheorem*{acka}{Acknowledgments}
\newtheorem{ack}{Acknowledgments}
\renewcommand{\theack}{}
\renewcommand{\appendixname}{\bf Appendix}
\renewcommand{\proof}{{\em Proof.\ }}

\hyphenation{
ap-pen-dix as-ymp-tot-ic at-trib-uted at-trib-ut-able
Bry-li-n-sky com-mu-ta-tion de-ge-ne-rate
de-riv-a-tive dis-trib-ute equi-vari-ant ex-tra-or-di-nary  
geo-met-ric griev-ance griev-ous grad-ed ho-lo-no-my ho-mo-thetic
in-fin-ite-ly in-fin-i-tes-i-mal Ha-rish Cha-n-dra mul-ti-plic-able 
non-euclid-ean non-iso-mor-phic non-smooth par-a-digm 
par-a-bol-ic pa-rab-o-loid pa-ram-e-trize phe-nom-e-non 
post-script pseu-do-dif-fer-en-tial pseu-do-fi-nite 
qua-drat-ics quad-ra-ture Han-kel rec-tan-gle semi-def-i-nite 
set-up wide-spread Euler-ian Feb-ru-ary Gauss-ian Grothen-dieck 
Hamil-ton-ian Her-mi-t-ian her-mi-t-ian Jan-u-ary 
Japan-ese Ka-shi-wa-ra Kor-te-weg Le-gendre No-vem-ber Rie-mann-ian 
Sep-tem-ber Za-mo-lo-d-chi-kov Kni-zh-nik quan-tum Op-dam
Mac-do-nald Ca-lo-ge-ro Su-ther-land Mo-ser 
Ol-sha-net-sky  Pe-re-lo-mov in-de-pen-dent ope-ra-tors 
cy-clo-to-mic ra-tio-nal de-gen-er-a-tion 
in-ter-est-ing de-for-ma-tions de-for-ma-tion pro-ce-dure 
fol-lows ope-ra-tors  pre-serve suf-fices ap-proach 
for-mu-las con-sider its com-ple-tion cor-re-spond-ing 
au-to-mor-phism be-cause pro-por-tional fi-nal-ly let-ting 
equi-v-a-lence ge-n-er-al-ized Mac-do-nald iden-ti-ties 
cor-re-s-pond sub-dia-grams par-ti-tion na-t-u-ral-ly 
or-dered stan-dard de-for-ma-tion ar-gu-ment com-bined 
sphe-r-i-cal rep-re-sen-ta-tions tri-go-no-me-t-ric
ge-n-er-al-ly speak-ing pri-m-it-ive ir-re-du-cible 
sum-ma-tion  rep-re-sen-ta-tives pro-por-ti-o-na-li-ty
ultra-sphe-ri-cal Ro-gers}

\def\ffor{\quad\hbox{ for }\quad}
\def\wwhen{\quad\hbox{ when }\quad}
\def\wwhere{\quad\hbox{ where }\quad}
\def\aand{\quad\hbox{ and }\quad}
\def\for{\  \hbox{ for } \ }
\def\iif{ \ \hbox{ if } \ }
\def\when{ \ \hbox{ when } \ }
\def\where{\  \hbox{ where } \ }
\def\and{\  \hbox{ and } \ }
\def\and{\  \hbox{ and } \ }
\def\oor{\  \hbox{ or } \ }
\def\proof{{\em Proof. \  }}

\def\equal{\stackrel{\,\mathbf{def}}{= \kern-3pt =}}

\def\la{\lambda}
\def\La{\Lambda}
\def\om{\omega}
\def\Om{\Omega}
\def\Th{\Theta}
\def\th{\theta}
\def\al{\alpha}
\def\be{\beta}
\def\ga{\gamma}
\def\ep{\epsilon}
\def\up{\upsilon}
\def\Up{\Upsilon}
\def\de{\delta}
\def\De{\Delta}
\def\ka{\kappa}
\def\kapp{\hbox{\bf \ae}}
\def\si{\sigma}
\def\Si{\Sigma}
\def\Ga{\Gamma}
\def\ze{\zeta}
\def\io{\iota}
\def\bio{b^\iota}
\def\aio{a^\iota}
\def\twio{\tilde{w}^\iota}
\def\hwio{\hat{w}^\iota}
\def\gio{\g^\iota}
\def\Bio{B^\iota}

\def\del{\delta}
\def\pa{\partial}
\def\vp{\varphi}
\def\ve{\varepsilon}
\def\inf{\infty}

\def\vph{\varphi}
\def\vps{\varpsi}
\def\vPh{\varPhi}
\def\vep{\varepsilon}
\def\vpi{{\varpi}}
\def\vth{{\vartheta}}
\def\vsi{{\varsigma}}
\def\vrh{{\varrho}}

\def\bph{\bar{\phi}}
\def\bsi{\bar{\si}}
\def\bvp{\bar{\varphi}}

\newcommand{\bS}{{\mathbf S}}
\newcommand{\bH}{{\mathbf H}}
\newcommand{\bF}{{\mathbf F}}
\newcommand{\bE}{{\mathbf E}}

\def\tal{\tilde{\alpha}}
\def\tbe{\tilde{\beta}}
\def\tde{\tilde{\delta}}
\def\tpi{\tilde{\pi}}
\def\txi{\tilde{\xi}}
\def\tPi{\tilde{\Pi}}
\def\tPhi{\tilde{\Phi}}
\def\tV{\tilde{V}}
\def\tJ{\tilde{J}}
\def\tla{\tilde{\lambda}}
\def\tga{\tilde{\gamma}}
\def\tGa{\tilde{\Gamma}}
\def\tvs{\tilde{{\varsigma}}}
\def\tu{\tilde{u}}
\def\tU{\tilde{U}}
\def\tw{\widetilde w}
\def\tW{\widetilde W}
\def\tB{\tilde B}
\def\tv{\tilde v}
\def\tV{\tilde V}
\def\tz{\tilde z}
\def\tb{\tilde b}
\def\ta{\tilde a}
\def\tih{\tilde h}
\def\trh{\tilde {\rho}}
\def\tx{\tilde x}
\def\tf{\tilde f}
\def\tg{\tilde g}
\def\tG{\tilde G}
\def\tk{\tilde k}
\def\tl{\tilde l}
\def\tL{\tilde L}
\def\tD{\tilde D}
\def\tR{\tilde R}
\def\tP{\tilde P}
\def\tH{\tilde H}
\def\tp{\tilde p}

\def\hH{\hat{H}}
\def\hh{\hat{h}}
\def\hR{\hat{R}}
\def\hY{\hat{Y}}
\def\hX{\hat{X}}
\def\hP{\hat{P}}
\def\hT{\hat{T}}
\def\hV{\hat{V}}
\def\hG{\hat{G}}
\def\hF{\hat{F}}
\def\hw{\widehat{w}}
\def\hW{\widehat{W}}
\def\hu{\hat{u}}
\def\hs{\hat{s}}
\def\hv{\hat{v}}
\def\hb{\hat{b}}
\def\hB{\widehat{B}}
\def\hze{\hat{\zeta}}
\def\hsi{\hat{\sigma}}
\def\hrh{\hat{\rho}}
\def\hth{\hat{\theta}}
\def\hy{\hat{y}}
\def\hx{\hat{x}}
\def\hz{\hat{z}}
\def\hg{\hat{g}}
\def\he{\hat{e}}
\def\hE{\widehat{E}}

\def\B{\mathbf{B}}
\def\I{\mathbf{I}}
\def\P{\mathbf{P}}
\def\G{\mathbf{G}}
\def\S{\mathbf{S}}
\def\F{\mathbf{F}}
\def\one{\mathbf{1}}
\def\Sn{\mathbf{S}_n}
\def\0{\mathbf{0}}
\def\H{\mathbf{H}}
\def\V{\mathbf{V}}

\def\f{\mathcal{F}}
\def\çF{\mathcal{F}}
\def\o{\mathcal{O}}
\def\t{\mathcal{T}}
\def\r{\mathcal{R}}
\def\l{\mathcal{L}}
\def\m{\mathcal{M}}
\def\k{\mathcal{K}}
\def\n{\mathcal{N}}
\def\d{\mathcal{D}}
\def\p{\mathcal{P}}
\def\cP{\mathcal{P}}
\def\a{\mathcal{A}}
\def\h{\mathcal{H}}
\def\c{\mathcal{C}}
\def\y{\mathcal{Y}}
\def\e{\mathcal{E}}
\def\v{\mathcal{V}}
\def\z{\mathcal{Z}}
\def\x{\mathcal{X}}
\def\s{\mathcal{S}}
\def\g{\mathcal{G}}
\def\u{\mathcal{U}}
\def\w{\mathcal{W}}
\def\i{\mathcal{I}}
\def\j{\mathcal{J}}
\def\b{\mathcal{B}}

\def\lan{\langle}
\def\llb{(\!(}
\def\ran{\rangle}
\def\rrb{)\!)}
 \def\dim{{\hbox{\rm dim}}_{\mathbb C}\,}
\def\lng{\hbox{\rm{\tiny lng}}}
\def\sht{\hbox{\rm{\tiny sht}}}
\def\sph{\hbox{\rm{\tiny sph}}}
\def\inv{\hbox{\rm{\tiny inv}}}

\def\br#1{\langle #1 \rangle}

\def\rank{\hbox{rank}}
\def\gl{\mathfrak{gl}_N}
%\def\sgn{\hbox{sgn}}
%\font\germ=eufb10 %at 12pt 
%\def\mathfrak#1{\hbox{\germ #1}}

\newcommand{\Aut}{\operatorname{Aut}}
\newcommand{\Hom}{\operatorname{Hom}}
\newcommand{\End}{\operatorname{End}}
\newcommand{\Ind}{\operatorname{Ind}}
\newcommand{\ad}{\operatorname{ad}}
\newcommand{\pr}{\operatorname{pr}}
\newcommand{\aweyl}{\tilde{\mathbb S}_n}
\newcommand{\hec}{{\mathcal H}^t_n}
\newcommand{\Func}{{\mathcal F}({\mathbb C}^n,{\mathcal H}^t_n)}
\newcommand{\tr}{\operatorname{tr}}
\newcommand{\Out}{\operatorname{Out}}
\newcommand{\Rad}{\operatorname{Rad}}
\newcommand{\Spec}{\operatorname{Spec}}
\newcommand{\id}{\operatorname{id}}
\newcommand{\Int}{\operatorname{Int}}
\newcommand{\ct} {\operatorname{ct}}

\newcommand{\rat}{{\mathbb Q}}
\newcommand{\real}{{\mathbb R}}
\newcommand{\cplx}{{\mathbb C}}
\newcommand{\zint}{{\mathbb Z}}

\newcommand{\sq}{\phantom{1}\hfill$\qed$}
\newcommand{\Rea}{\Re}
\newcommand{\Ima}{\Im}

\newcommand{\st}{\bowtie}
\newcommand{\modd}{\mbox{\,mod\,}}
\newcommand{\lr}{\langle}
\newcommand{\rr}{\rangle}
\newcommand{\eps}{\varepsilon}
\newcommand{\phk}{\phi^{(k)}}
\newcommand{\psk}{\psi^{(k)}}
\newcommand{\Res}{\mbox{Res}\;}
\newcommand{\sgn}{\mbox{sgn}}
\newcommand{\mn} {\left\{ \begin{array}{c}m\\
n\end{array}\right\}}

\def\sX{\mathscr{X}}
\def\sH{\mathscr{H}}
\def\sY{\mathscr{Y}}
\def\TT{\mathfrak{T}}
\def\JJ{\mathfrak{J}}
\def\HH{\mathfrak{H}}
\def\FF{\mathfrak{F}}
\def\GG{\mathfrak{G}}
\def\CC{\mathfrak{C}}
\def\LL{\mathfrak{L}}

\def\BB{\mathfrak{B}}
\def\AA{\mathfrak{A}}
\def\ZZ{\mathfrak{Z}}
\def\HH{\hbox{${\mathcal H}$\kern-5.2pt${\mathcal H}$}}
\def\HHH{\hbox{${\mathbb H}$\kern-4.2pt${\mathbb H}$}}
\def\tHH{\widetilde{\HH\ }}

\font\smm=msbm10 at 12pt 
\def\symbol#1{\hbox{\smm #1}}
\def\lsmash{{\symbol n}}
\def\rsmash{{\symbol o}}
\def\#{\sharp}

\font\tenbf=cmbx10
\font\tenrm=cmr10
\font\tenit=cmti10
\font\ninebf=cmbx9
\font\ninerm=cmr9
\font\nineit=cmti9
\font\eightbf=cmbx8
\font\eightrm=cmr8
\font\eightit=cmti8
\font\sevenrm=cmr7
\font\sevenbf=cmbx7

%END MACROS

\title [Punctual Hilbert schemes in dimension two]
{Gr\"obner cells of punctual Hilbert schemes in dimension two}
\author[Ivan Cherednik]{Ivan Cherednik $^\dag$}
%\date{February 2, 2014}

\begin{abstract}
We begin with a comprehensive discussion
of  the punctual Hilbert scheme
of the regular two-dimensional local ring 
in terms of the Gr\"obner cells. These schemes are 
the most degenerate fibers of the Grothendieck-Deligne norm map
(the Hilbert-Chow morphism), playing an important role in the 
study of Hilbert schemes of smooth surfaces.
% (G\"otsche, Haiman, Nakajima and others). 
They are generally singular, but their Gr\"obner cells 
are affine spaces; they 
admit an explicit parametrization due to Conca and Valla. 
We use this to obtain the  Gr\"obner decomposition
of compactified Jacobians of plane 
curve singularities, which is non-trivial even for
the generalized Jacobians (principal ideals only).
One of the application is the topological invariance of
certain variants of compactified Jacobians and 
the corresponding motivic superpolynomials for analytic 
deformations 
of quasi-homogenous plane curve singularities and some
similar families.
%Another 
%application, is to the wall-crossing based on weighted 
%Gr\"obner schemes.  
\end{abstract}

\thanks{$^\dag$ \today.
\ \ \ Partially supported by NSF grant
DMS--1901796 and ERC grant AdG 669655.  }

\address[I. Cherednik]{Department of Mathematics, UNC
Chapel Hill, North Carolina 27599, USA\\
chered@email.unc.edu}

 \def\sht{\raisebox{0.4ex}{\hbox{\rm{\tiny sht}}}}
 \def\bysame{{\bf --- }}
 \def\~{{\bf --}}

%%%rk>1 PAPER ONLY:
\def\rk{r\!k} % \def\rk{{\mathsf \varrho}}
\renewcommand{\dim}{\hbox{dim\,}}   %%%OTHERWISE dim_C

 \def\rr{{\mathsf r}}
 \def\cc{{\mathsf c}}
 \def\ss{{\mathsf s}}
 \def\mm{{\mathsf m}}
 \def\pp{{\mathsf p}}
 \def\ll{{\mathsf l}}
 \def\aa{{\mathsf a}}
 \def\bb{{\mathsf b}}
 \def\NS{\hbox{\tiny\sf ns}}
 \def\ssum{\hbox{\small$\sum$}}
\newcommand{\comment}[1]{}
\renewcommand{\tilde}{\widetilde}
\renewcommand{\hat}{\widehat}
\renewcommand{\V}{\mathbb{V}}
\renewcommand{\F}{\mathbb{F}}
\newcommand{\dagx}{\hbox{\tiny\mathversion{bold}$\dag$}}
\newcommand{\ddagx}{\hbox{\tiny\mathversion{bold}$\ddag$}}
\newtheorem{conjecture}[theorem]{Conjecture}
\newcommand*\toeq{
\raisebox{-0.15 em}{\,\ensuremath{
\xrightarrow{\raisebox{-0.3 em}{\ensuremath{\sim}}}}\,}
}
\newcommand{\unknot}{\hbox{\tiny\!\raisebox{0.2 em}{$\bigcirc$}}}
%\WarningFilter{latex}{Text page} %%%MAKE IT % FOR ARXIV
\newcommand*{\vect}[1]{\overrightarrow{\mkern0mu#1}}
\newcommand*{\medcap}{\mathbin{\scalebox{.75}
{\ensuremath{\bigcap}}}}
\newcommand*{\medcup}{\mathbin{\scalebox{.75}
{\ensuremath{\bigcup}}}}
%\WarningFilter{latex}{Text page} %%%MAKE IT % FOR ARXIV

\vskip -0.0cm
%\par
%{\centering
%\medskip
%\par}
%\vskip -0.0cm
\maketitle
\vskip -0.0cm
%\smallskip

\vskip 0.2cm
\noindent
{\em\small {\bf Key words}: Hilbert schemes, affine plane,
Grothendieck-Deligne map, \\
Gr\"obner cells, zeta functions, plane curve 
singularities.}
\smallskip

{\tiny
\centerline{{\bf MSC} (2010): 14H50, 17B22, 17B45, 20C08, 20F36,
22E50, 22E57, 30F10, 33D52, 33D80, 57M25.}
}
\smallskip

\vskip -0.5cm
\renewcommand{\baselinestretch}{0.95}
{\small
%%%\pagenumbering{gobble}
\tableofcontents   %%% REMOVE % IN ARXIVE!!!!
}
\renewcommand{\baselinestretch}{1.0}

\renewcommand{\natural}{\wr}

%\vfill\eject
\setcounter{section}{0}
\setcounter{equation}{0}
\section{\sc Introduction}
We begin with a discussion of the Hilbert scheme $H^{(n)}$ of 
the local ring $\C[[x,y]]$ in terms of Gr\"obner cells $C_\la$, 
providing all details.  
Here $n$ is the codimension (the length) of the ideals in 
$\C[[x,y]]$ and $\la$ is a partition of $n$. 
The cellular decomposition of $H^{(n)}$ was obtained
\cite{ES} via \cite{BB} 
from the analysis of the action of the maximal torus in 
$SL(3,\C)$ in the tangent 
spaces of $H\!ilb^{(n)}(P^2_{\C})$ at the corresponding
fixed points, which are monomial ideals. Our starting point is an
entirely local definition of the Gr\"obner cells 
of $H^{(n)}$ and their explicit parametrization following
\cite{CV}, Theorem 3.3 ($i=2$).
%For the convenience of the
%readers, $H\!ilb^{(n)}$ is connected with 
%$H\!ilb^{(n)}(\C^2)$ in full detail. 
The Gr\"obner decomposition is an important tool in the study
of $H^{(n)}$; it can be potentially
used for any isolated surface singularities, not only
quasi-homogeneous. 
%\vfil

We note that the definition of the Gr\"obner cells generally 
requires some {\it valuation\,} of  $\C[x,y]$ or $\C[[x,y]]$. 
The dependence on the ratio of the valuations
of $x$ and $y$ can be interpreted as some
{\it wall-crossing\,}.

This decomposition induces
that of compactified Jacobians of plane 
curve singularities under the natural embeddings,
which is  non-trivial even for
the generalized Jacobians, which are for
{\it principal\,} fractional
ideals instead of all of them.
Combined with that in terms of the
so-called Piontkowski strata, it can be presumably used
to obtain the super-duality for the motivic superpolynomials 
from \cite{ChP1}. 
This duality was conjectured in \cite{Ch} to coincide with
the functional equation of the Galkin-St\"ohr zeta functions
for any plane curve singularities; this is related to 
motivic theory of $H^{(n)}$.
See \cite{Sto} and also  \cite{MY,MS,ORS}. 
%\vfil

The topological invariance
of these superpolynomials is a significant part of
the conjectures in \cite{ChP1,Ch}. We prove a stronger
fact for the deformations of quasi-homogeneous singularities
and some similar families. Namely, 
$\tilde{Jac^\bullet}$ from Lemma \ref{lem:tilde-Jac},
closely related to the compactified Jacobians,
are topological invariants for such families.

\subsection{Punctual Hilbert schemes}
The Hilbert scheme $H^{(n)}$ of 
$\C[[x,y]]$ is defined as a scheme of all its ideals
of codimension $n$. Equivalently, it can be introduced
as the fiber of the  Grothendieck-Deligne 
norm map (the Hilbert-Chow morphism) 
$\pi_n:\, H\!ilb^{(n)}(\C^2)\to S^n(\C^2)$  
over the point $nO$ in the symmetric power  $S^n(\C^2)$ of the
affine plane $\C^2$; see \cite{Del}.
Here $O=\{x=0,y=0\}$, $H\!ilb^{(n)}(\C^2)$ is
the Hilbert scheme formed by ideals $I\subset \C[x,y]$ of
codimension $n$. This fiber is called the {\it punctual Hilbert 
scheme}; some authors call them "local punctual".
\vskip 0.2cm

Let $(x,y)_\ep=(\ep^u x,\ep^v y)$
be the action of the torus $\C^*\ni \eta $ in $\C^2$, depending
on $u,v\in \R_+$. Assuming that $0<u<<v$ and $\ep>0$, 
the limit  $I^0=\lim_{\ep\to 0} I_h$ is well-defined and is 
a {\em monomial ideal\,} in  the same $H\!ilb^{(n)}(\C^2)$. 
Monomial ideals are those linearly generated by $x^a y^b$. 
Combinatorially, we can obtain $I^0$ simply by taking from any 
$f\in I$ its top monomial $f^0$ under the
following lexicographic ordering:
$1<y<y^2<\cdots <x<xy<xy^2<\cdots$.  
\vfil

The monomial ideals 
are fully determined by the partitions 
$\la=\{m_1\ge m_2\ge \cdots \ge m_\ell>0\}$, where 
$\sum_{i=1}^\ell m_i=n$. Let $I_\la$ be linearly generated
by $x^a y^b$ such that $\{a,b\}\in \Z_+^2\setminus \la'$, 
where we represent $\la$ as the set 
$\la'\equal
\bigl\{\{i,j\}\in \Z_+^2 \mid 0\le i <\ell,\, 
0\le j < m_{i+1}\bigr\}$.
 Equivalently, 
$x^a y^b\in I_\la$ if and only if  
$a\ge i^{\circ}$ and $b\ge j^{\circ}$ for at least one {\em corner\,}
$\{i^{\circ}, j^{\circ}\}$ of $\Z_+^2\setminus \la'$. 
Any monomial ideal is $I_\la$ for some $\la$.
\vfil

For a given partition $\la\vdash n$, the {\em Gr\"obner cell\,} 
$Gr_\la$ of $H\!ilb^{(n)}(\C^2)$  is defined as follows: 
$Gr_\la\equal \{I\in H\!ilb^{(n)}(\C^2) \mid I^0=I_\la\}.$
They form a {\it cellular
decomposition\,} of $H\!ilb^{(n)}(\C^2)$ in the sense of Fulton.
\vfil

Using the embedding $H^{(n)}\subset H\!ilb^{(n)}(\C^2)$, let
 $Gr^0_\la\equal Gr_\la \cap H^{(n)}$. There is an entirely
local definition of the Gr\"obner cells for $H^{(n)}$.
Namely, we switch the order of $x$ and $y$ (now $y>x$), and take 
the {\em lowest monomials\,} instead of the top ones used
in the construction of $I^0$.
The notation is $I_0$ instead of $I^0$ throughout the paper.
The corresponding strata of $H^{(n)}$ will be  denoted by
$C_\la$. This definition is compatible with the passage to the
completion $\C[[x,y]]$ of $\C[x,y]$ at $(x=0=y)$, i.e.
it is indeed local.  
\vfil

There is an embeddings $C_\la\hookrightarrow Gr_\la$,
which results in the identification $C_\la\simeq Gr_\la^0$.
We make this very explicit. Note that
the cell decomposition of $H^{(n)}$ is obtained in \cite{ES}
without any reference to $H\!ilb^{(n)}(\C^2)$. It is some
part of the decomposition of $H\!ilb^{(n)}(P_{\C}^2)$, and there
is some Poincar\'e duality between $H^{(n)}$ and
$H\!ilb^{(n)}(\C^2)$ (see below).

\subsection{Some basic facts}
In contrast to  
$H\!ilb^{(n)}(\C^2)$, the scheme $H^{(n)}=\pi_n^{-1}(nO)$ is projective
and generally not smooth. It is irreducible due to J.~Briancon 
\cite{Bri} and of
dimension $n-1$; see also \cite{Ia}. Furthermore,
it is a complete intersection and reduced, so it is
Cohen-Macaulay. This is due to M.~Haiman; see Proposition 2.10 
in \cite{Ha1}.

Here $\C^2$ can be replaced by any
{\em smooth\,} quasi-projective surface $X$; 
the corresponding  $H\!ilb^{(n)}(X)$ formed by    
subschemes in $X$ supported in one (any) point is
isomorphic to $H^{(n)}$, i.e. to that for $\C^2$. 
The fibers $\pi_n^{-1}(nP)$  are the most degenerate
ones; knowing them
is sufficient to calculate {\em any\,} fibers of $\pi_n$
for any $X$. 

Namely, 
let $\mu=\{n_1\ge n_2\ge \cdots \ge n_r>0\}$ be a partition of $n$,
$P_1,\ldots, P_r$ be pairwise distinct points. We set
$D_\mu=\sum_{i=1}^r n_iP_i$. Then $\pi_n^{-1}(D_\mu)$ is naturally
isomorphic to the product of $\pi_{n_i}^{-1}(n_iO)$ over 
$1\le i\le r$.  We arrive at the fibration 
of $H\!ilb^{(n)}(X)$ with respect to the standard 
stratification of $S^{n}(X)$ with the strata 
$\d_\mu\equal \{D_\mu\}$. The latter are unramified covers of
$\prod S^{r}(X)$ minus the diagonals, i.e. they are
smooth. The fibration of $H\!ilb^{(n)}(X)$ corresponding to
$\pi_r$ is locally trivial upon the restriction to any
$\d_\mu$.
\vskip 0.1cm
\vfil

From $\dim H^{(n)}\!=n\!-\!1$, we obtain
that $\dim \pi_n^{-1}(D_\mu)=
 \sum_{i=1}^r (n_i-1)=n-r$, i.e. it is $|\mu|$ minus the 
length of $\mu$. This formula can be possibly
connected with the formula $\dim C_{\la}=n-\ell(\la)$ (below)
via some deformation procedure.  
The fibers $H^{(n)}$ are the key in the theory of 
$H\!ilb^{(n)}(\C^2)$ and its various applications; 
see e.g. \cite{Ha2}. The explicit parametrization
of $C_\la$, can be helpful here. 
\vskip 0.2cm
\vfil

Importantly, we do not
need the action of $\C^*$ in the definition of $C_\la$,
and it is  entirely local. It can be extended to more general
isolated surface singularity, not only quasi-homogeneous.
An explicit parametrization of the corresponding cells can be
involved, but the direct calculations follow the same lines as
for $\C[[x,y]]$ 
and are doable for relatively simple surface singularities. 
\vskip 0.2cm

Using the smoothness of $H\!ilb^{(n)}(\C^2)$ (J.~Fogarty)
and considering
its tangent space at $I_\la$, the general result of 
A.~Bialynicki-Birula \cite{BB} gives that $Gr_\la$ is 
an affine space. Here the definition of $I^0$ 
via the action of $\C^*$ is used; this action can be calculated
explicitly in the tangent space at $I_\la$. For instance,
it gives that
the dimension of $Gr_\la$ is $n+m(\la)$, 
where $m(\la)=m_1$. See \cite{ES}, Theorem 1.1, ($iii$),
\cite{CV}, Theorem 3.3 ($i=2$), 
and also  \cite{Nak1,Nak2,MO} for different aspects and
generalizations. 
Note that the embedding 
$C_\la \hookrightarrow Gr_\la$
is from the space of dimension $n-\ell(\la)$ to the one
of dimension $n+m(\la)$.

\vskip 0.2cm
Generally, \cite{BB} cannot be used for $H^{(n)}$ since it 
is not smooth. However, G.~Ellingsrud and S.~Str{\oin}mme
obtain a cellular decomposition of $H^{(n)}$ 
as part of that of
$H\!ilb^{(n)}(P_{\C}^2)$, which is smooth. 
The cells are affine spaces, which follows from \cite{BB}.
They calculate Betti numbers
of $H^{(n)}$  in their 
Theorem 1.1, ($iv$), which gives the number of cells
and their dimensions. We note that
taking here the $n$-row (in our notations) as  $\la$ readily gives the 
irreducibility of $H^{(n)}$; see Corollary 1.2 in \cite{ES}.
Namely, the closure of the corresponding  \Yboxdim5pt
$C_{\yng(2)\cdots}$ of dimension $n-1$, 
the {\em big cell\,}, is the whole $H^{(n)}$.

\subsection{Hilbert-type zetas} 
The classical {\em Hasse-Weil zeta} has the following
presentation: 
$Z(X,t)=\sum_{n=0}^\infty t^n |S^n(X)(\F_q)|$.
This formula holds for any varieties
$X$ over $\F_q$ and their symmetric powers $S^n(X)$.
Here $X$ can be singular; then $0$-cycles over $\F_q$ must 
be counted instead of the points. See e.g. \cite{Mus}.

Importantly, the Betti numbers 
$b_{i}(X)$, the ranks of Borel-Moore $i${\tiny th} homology,
are related to the classical zeta-function for
smooth $X$. Algebraically, they occur
as the degrees of the polynomials in $t$
for the contributions of the corresponding cohomology 
in Weil's formula.

The knowledge of the cell decomposition of $H^{(n)}$ generally
can be used to calculate the {\em Hilbert-type\,} zeta-functions of 
{\em any\,} quasi-projective smooth surfaces $X$, which are defined
as follows. We replace $S^n(X)$ by $H\!ilb^{(n)}$, setting
$\z(X,t)=\sum_{n=0}^\infty t^n |H\!ilb^{(n)}(X)(\F_q)|$.
In its {\em motivic\,} counterpart, the cardinality
$|H\!ilb^{(n)}(X)(\F_q)|$
is replaced by the class of $H\!ilb^{(n)}(X)$  in
the Grothendieck ring of varieties over
the basic field; the count of $\F_q$-points is then considered 
as the {\em counting motivic measure}. 
\vskip 0.2cm

Let us outline the general way of calculating
such $\z$-functions based on Weil conjectures.
The {\em Hasse-Weil zeta} $Z(X,t)$ gives
all $|S^n(X)(\F_q)|$, which is sufficient to  
calculate all $|\d_\mu(\F_q)|$ using the inclusion-exclusion
principle. Then we can apply the 
formulas for $|H^{(n)}(\F_q)|$: each cell $\mathbb{A}^m$
of $X$ results in $q^m$. This gives $\z(X,t)$. 
We will provide 
an example of such a calculation below.
This is related to the way L.~G{\"o}ttsche obtained
his well-known formula from \cite{Got1,Got2} 
in terms of Betti numbers of 
Hilbert schemes of smooth quasi-projective surfaces $X$. 
\vskip 0.2cm

A significant part of the paper is devoted to
the connections with the plane curve singularities.
%aimed
%at their {\em motivic superpolynomials\,} from \cite{ChP1}.
For any element $P(x,y)$ of $I\in H^{(n)}$, this ideal 
is the inverse image of some ideal in $\r=\C[[x,y]]/(P(x,y))$
of codimension $n$. Any such $I$  
has a {\em canonical\,}
"first Gr\"obner generator" $P_I(x,y)$, the one with the minimal
possible $x^a (a>0)$ and other monomials in it
"from" the boxes of  $\la$. We fix $P(x,y)$ and identify
the subset $\{I\in C_\la\,\mid\, P_I(x,y)=P(x,y)\}$
with some subset of the compactified Jacobian 
$\overline{Jac}_{\r}$ of $\r$. Cf. Theorem \ref{thmtoI}, $(ii)$. 
%the Gr\"obner cells for $\C[[x,y]]$ with Hilbert schemes of $\c$
%and their compactified Jacobians $\overline{Jac}_{\,\c}$. 
Using the parametrization of $C_\la$, this subset is 
given by explicit equations. 

We call such
subsets {\em Gr\"obner
strata} of $\overline{Jac}_{\r}$, which stratification
can be related to \cite{MY,MS}.  The strata can be singular
as schemes, 
but they are affine spaces topologically in many cases;
see  (\ref{C-incl}) and around,
This decomposition can be 
quite non-trivial even in the case of the {\em generalized Jacobian},
$Jac_{\r}\subset \overline{Jac}_{\r}$, 
the group of invertible $\r$-submodules in the normalization
ring of $\r$. 
% At least in examples,
%there is a connection with the super-duality of 
%$\z$-functions and {\em motivic superpolynomials} of
%plane curve singularities; see e.g. \cite{Ch}. 
\vskip 0.2cm

 \Yboxdim5pt
{\sf Example}. Let us show how to calculate 
$|H\!ilb^{(3)}(\F_q)|$ for $X=\mathbb A^2$ using this approach. 
We need the formulas $|Gr_\la(\F_q)|=q^{n+m(\la)},
|C_\la(\F_q)|=q^{n-\ell(\la)}$, and the number of
$\F_q$\~points of the fiber in $H\!ilb^{(n)}$
over $n_1P_1+\cdots+ n_rP_r\in S^n(X)$ over $\F_q$ for
$\mu=\{n_1\ge n_2\ge \cdots \ge n_r>0\}\vdash n$ and
pairwise distinct points $P_1,\ldots, P_r$. The latter equals\ 
$\prod_{i=1}^r |H^{(n_i)}(\F_q)|$.

We already know that $|H\!ilb^{(3)}(\F_q)|=\sum_{\la\,\vdash 3}
|Gr_\la(\F_q)|=q^6+q^5+q^4$. Let us obtain this quantity
using the approach via $S^n(X)$, which is generally applicable
to any smooth surfaces $X$ over $\F_q$. Using
the formula $|H^{(m)}(\F_q)|=\sum_{\la\,\vdash m}
|C_\la(\F_q)|$, we obtain:  $|H^{(1)}(\F_q)|=1$,
 $|H^{(2)}(\F_q)|=1+q$, $|H^{(3)}(\F_q)|=1+q+q^2$.
\vskip 0.3cm

Finally, showing the source of the terms in 
$\lan\!\lan\cdots\ran\!\ran$:
\begin{align*}
&|H\!ilb^{(3)}(\F_q)|=(1+q+q^2)q^2\, 
\lan\!\lan\mu=\yng(3)\,\ran\!\ran +
(1+q)q^2(q^2-1)\, \lan\!\lan \mu=\yng(2,1)\,\ran\!\ran\\
&+\frac{q^2(q^2-1)(q^2-2)}{6}\,
\lan\!\lan\mu=\yng(1,1,1)\,,  P_i\in \F_q \ran\!\ran +  
q^2\frac{(q^4-q^2)}{2}\\
&\lan\!\lan P_1\in \F_q, P_{2,3}\in\F_{q^2}\setminus \F_q\ran\!\ran+
\frac{(q^6-q^2)}{3}\, \lan\!\lan P_i\in \F_{q^3}\setminus 
\F_q \ran\!\ran=
q^6+q^5+q^4.
\end{align*}
\vskip 0.2cm

\subsection{Toward functional equation}
Inspired by the theory of plane curve singularities, it was
expected in \cite{Ch} that 
$\l$-functions of reasonably good isolated surface singularities 
over $\F_q$ depend on $q$ uniformly, which property is called 
"strong polynomial count", and satisfy the functional equation. 
These $\l$\~functions are infinite products in contrast to 
those for plane curve singularities, so the functional equation
will be with "infinite" scaling factors. We will provide
examples below. Here only Hilbert-type zeta-functions make
sense to consider; symmetric powers 
of a singularity (just a point) are meaningless. 

An expected 
connection with the $q$-deformations of the classical $L$-functions 
from Number Theory is touched upon in \cite{Ch}; this is quite a
motivation, but very preliminary by now.  
\vskip 0.2cm

The definition of the {\it local} zeta-function of
any singularity ring $\r$
is straightforward: $\z_{\r}(t)\equal 
\sum_{n=0}^\infty t^n |H_{\r}^{(n)}(\F_q)|$, where
$H_{\r}^{(n)}$  is a scheme of ideals in $\r$ of codimension $n$.
and $\r$ must be at least Gorenstein. The latter
can be insufficient: it is  expected that
the surface isolated singularities corresponding to  
Seifert $3$-folds (as their links) constitute a natural class. 
The uniform dependence on $q$ is an important test, what is
called "strong polynomial $q$-growth" of $|H_{\r}^{(n)}(\F_q)|$.
\vfil

The rings $\r$ are initially over $\C$, so we need to consider 
them over proper extensions of $\Z$, and then switch to $\F_q$,
where $q=p^m$ assuming that primes $p$  are of
"good reduction"; almost all $p$ are such. This passage to
$\F_q$ is sufficiently well understood for curve singularities. 
Then we switch to the $\l$-functions; the functional 
equation is expected only for them. 
\vskip 0.2cm
\vfil

{\sf Plane curve singularities.}
For the rings $\r$ of {\em Gorenstein 
curve singularities}, one has:
 $\l_{\r}\equal (1-t)\z_{\r}(t)$,
which is the Galkin-St\"ohr $L$-function. 
It is a polynomial in terms of $t$; see \cite{Sto}. 
The multiplication by $(1-t)$ is a counterpart of
the multiplication by $(1-t)(1-qt)$ for smooth projective
curves.  

Conjecture 4.5 from \cite{Ch}
states that $\l_{\r}=\h_{\r}(qt,t)$ for
{\em plane curve singularities\,} $\r$, where
the {\em motivic superpolynomial} $\h_{\r}(q,t)$
is as follows. We consider $\r$ as a subring of
$\F_q[[z]]$, where $z$ is the uniformization
parameter; so $\r$ is an arbitrary subring in $\o$
with two generators and such that its fields of
rationals coincide for $\r$ and $\F_q[[z]]$. Then 
$\h_{\r}(q,t)\equal
\sum_{M}t^{dim_{\F_q}(\F_q[[z]]/M)}$, where the
summation is over $\r$-submodules $M\subset \F_q[[z]]$ such that
$M\F_q[[z]]=\F_q[[z]]$. 

Generally,  the coincidence $\l_{\r}=\h_{\r}(qt,t)$
does not hold for non-planar Gorenstein
curve singularities, so the usage of Punctual 
Hilbert schemes $H^{(n)}$ and similar
objects seems inevitable here. 

The substitution $q\mapsto qt$ requires an
assumption that
$\h_{\r}(q,t)$ is a polynomial in terms of $q$, 
conjectured for any plane curve singularities.  
Then a conjectural relation with the DAHA superpolynomials
from \cite{ChP1} gives
that $\h_{\r}$ are
{\em topological invariants\,}, i.e. depend only on the
corresponding  valuation semigroup $\Ga_\r$; see Section \ref{sec:CompJ}.
This will follow from our considerations
for some families of $\r$.
%\vfil
\vskip 0.2cm
{\sf Surface singularities.}
For surface singularities, 
the division of  $\z_{\r}(t)$ by
$\z_{\o}(t)$ for $\o=\C[[x,y]]$ can be expected, but this 
can be more involved than this.  Anyway, the objective is to have the
functional equation for $\l_\r(t)$
with respect to the substitution $t\mapsto 1/(q^2t)$.

To explain which functional equation can be expected here,
let us reproduce the G\"ottsche formula:
$$
\sum_{n\ge 0}
\sum_{i\ge 0} (-1)^n b_i(H\!ilb^{(n)}(X))q^{i/2} t^n=
\prod_{k\ge 1}\prod_{j=0}^4
(1-q^{k-1+j/2}t^k)^{(-1)^{j+1}b_j(X)}
$$
for the Betti numbers of $X$ and $H\!ilb^{(n)}(X)$. 
Using the Poincar\'e duality for smooth projective
$X$:\  $b_j(X)=b_{4-j}(X)$ and 
the right-hand side satisfy
a {\em formal\,} functional equation
upon the substitution $t\mapsto 1/(q^2t)$ from
the classical functional equation for surfaces. This of course
holds up to some {\em infinite (!)\,} monomial in terms 
of $q^{1/2}$ and $t$, necessary to get rid of the 
denominators in the binomials. Such an "infinite rescaling"    
will not be addressed in this paper. 
  
Here $q$ is treated as a free parameter.
Let us assume that all $H\!ilb^{(n)}(X)$ have {\em cellular 
decompositions\,}; see e.g. Proposition 1.5 from \cite{ES}.
Then the left-hand side above coincides with $\z(X,t)$, based on
counting the $\F_q$-points.  Without this assumption, we have 
generally two different approaches, the modular one and its geometric
counterpart based on the Borel-Moore homology
(homology with closed support) or other kinds of
(co)homology. The same functional equation is expected for
either one, however the tools for its verification will be 
very different. Both approaches can be potentially used for  
$\l$-functions of (reasonably good) isolated 
surface singularities.
%\vskip 0.2cm
\vfil

{\sf Back to $\C^2$.} For $X=\mathbb{A}^2$ considered over $\F_q$,
 one has:
$$
\z(X,t)=\!\sum_{n=0}^\infty t^n |H\!ilb^{(n)}(X)(\F_q)|
=\!\sum_{n=0}^\infty\! \sum_{\,\la\,\vdash n} q^{n+m(\la)}=
\prod_{i=1}^\infty (1-q^{i+1}t^i)^{-1}.
$$
See \cite{ES} and Remark 4.7 in \cite{KR}. 
The formula for the local Hilbert-type
zeta of $\r=\F_q[[x,y]]$ is
$$
\z_\r(t)=\!\sum_{n=0}^\infty t^n |H^{(n)}(\F_q)|
=\!\sum_{n=0}^\infty\!\sum_{\,\la\,\vdash n} q^{n-\ell(\la)}=
\prod_{i=1}^\infty (1-q^{i-1}t^i)^{-1}.
$$

The similarity of the latter with that for $\mathbb{A}^2$
is not accidental. Following 
 \cite{ES},  let us use the decomposition
$\mathbb{P}^2=\mathbb{A}^2\cup
\mathbb{A}^1\cup \mathbb{A}^0$. Accordingly,
$\z(\mathbb{P}^2,t)$ is the product of the corresponding
zetas for the $0$-dimensional subschemes in $\mathbb{P}^2$ 
supported in $\mathbb{A}^2, \mathbb{A}^1, \mathbb{A}^0$. The 
products above are for $\mathbb{A}^2$ and $\mathbb{A}^0$.
They can be readily seen in the G\"ottsche formula
for $\mathbb{P}^2$: the products corresponding to $b_0$ and $b_4$. 
Thus they are dual to each other with respect to the 
Poincar\'e duality, which is here the 
{\em formal\,} substitution
$t\mapsto 1/(q^2t)$  followed by the rescaling, the
multiplication by
an {\em ininite\,} $q,t$\~monomial. 
 The product  for $\mathbb{A}^1$ is in
terms of $q^it^i$; see \cite{ES}. It corresponds to
$b_{2}(X)$; all odd Betti numbers are $0$.

\subsection{Some perspectives}
Let us mention here a series of papers on the generating 
function for
Euler numbers of Hilbert schemes of points on
a simple singularities  $\C^2/\Ga$ for finite subgroups
$\Ga\subset SL(2,\C)$. The conjecture 
by {A.~Gyenge}, {A.~N{\'e}methi}, and {B.~Szendr{\"o}i}, 
was that it is the
character of the corresponding Kac-Moody 
basic representation, where its torus variable
is evaluated at a proper root of unity of order related
to the Coxeter number.  It was checked  in \cite{DS,Tod}
in type $A$,  and for $D,E$ in \cite{GNS1,GNS2,Nak3}. 
See Theorem 1 from the last reference for an exact statement. 
This is directly related to our generating function,
but we need the {\em refined\,} version of this theorem
in terms of the Betti numbers. Also, we focus on the
local theory of isolated surface singularities, i.e. on 
{\em punctual\,} Hilbert schemes. 

Another equally important direction is the passage to
{\em instantons}, i.e. to torsion-free sheaves of
any ranks instead of ideals in the definition of Hilbert
schemes. Here we have a solid theory in any ranks
for plane curve singularities from 
\cite{ChP2}, though with quite a few conjectures.
Even the cases of affine plane and its local version 
for $\r=\C[[x,y]]$ are quite interesting, directly
related to the {\em Nekrasov instanton sums}. More generally, 
such theory must be associated with arbitrary Young diagrams;
the sheaves of rank $r$ correspond to the $r$-column. The
functional equation then includes the transposition of the Young
diagrams. 

The generalization to {\em arbitrary\,} Young diagrams 
is done by now for a {\em different} DAHA-based approach,
which conjecturally gives the same superpolynomials as those
from the motivic theory in \cite{ChP1,ChP2} (for columns). Also, 
let us at least mention the
connections of the  compactified
Jacobians to {\em affine Springer fibers}; 
see e.g. \cite{Yu}. For this, we make $a=0,t=1$ in the superpolynomials,
so the functional equation, which
requires $t$, generally can not be seen. The parameter $a$, which
we do not introduce in this paper, is associated with complete
{\em flags} of the modules and the ideals, related to the {\em nested
Hilbert schemes}. Importantly, $a\mapsto a$ in the functional equation. 
%\vfil

\vskip 0.2cm 
The totally local theory has many advantages;
torsion-free bundles over singular curves and surfaces are
a difficult topic in classical algebraic geometry. In the case of
local isolated singularities, the corresponding definitions
are actually no different from those in \cite{ChP2} 
in dimension one. Presumably the local approach captures many
features of the theory of instanton sums and related directions.
It is important that the  
$\z$-function  of $\C[[x,y]]$ in terms of $\{H^{(n)}\}$,
is closely related  to that of $\C^2$; it is certainly no simpler
in spite of its local nature.  
%\vfil

\vskip 0.2cm
A related direction is the {\em combinatorial
wall-crossing}.
One can generalize the Gr\"obner decomposition, coupling it
with the valuations of $\C[[x,y]]$. For instance, let 
$val_{r,s}(x^a y^b)=ra+sb$ for relatively prime numbers 
$r,s\in \Z_+$, $val_{r,s}(p)$ be the minimum
of valuations  of monomials in a polynomial $p(x,y)$.
Then we do the following. First, we group and order
the monomials in any 
$f\in \C[[x,y]]$ with respect to $val_{r,s}$. Second,
we find the minimal monomial in the groups with coinciding
valuations  for 
the Gr\"obner ordering $\{x^\infty\!<\!y\}$. 

For $r=0,s=1$, this gives the standard definition; 
the case $r=1,s=0$ corresponds to the switch 
of $x$ and $y$. Given a valuation, the corresponding Gr\"obner 
strata provide some basic elements in Borel-Moore homology 
of $H^{(n)}$. The construction depends on $r/s$, so we have 
{\it connection matrices}. We will not develop this   
in the present paper, but some relations to 
the singularities $\{x^s=y^r\}$ will be discussed.

%Compare with \cite{EH} and see the references there.

Generally, allowing the valuations to be $\infty$, we naturally
arrive at plane curve singularities; namely,
$\r=\C[[x,y]]/\{f\in \C[[x,y]]\mid  val(f)=\infty\}$.
The semigroup of all valuation gives its topological
type in the unibranch case. In a sense, {\em tropical geometry}
is when $0$ is allowed.

All these and related directions obviously require as 
constructive theory of $H^{(n)}$ as possible, which is the 
subject of the present paper.

\smallskip
{\bf Acknowledgements.}
The author is very thankful to David Kazhdan and Mikhail Finkelberg
for useful discussions.
Special thanks to Giovanni Felder for valuable contributions
to the paper concerning the Gr\"obner cells of $\C^2$
and related matters. The referee's comments greatly
helped to improve the paper and fix some problems, 
especially in Theorem \ref{thmc} and around
Proposition \ref{MIast}.
And many thanks to Rahul Pandharipande and
ETH-ITS (Zurich) for hospitality.

\setcounter{equation}{0}
\section{\sc Gr\"obner cells}
We will begin with some connections between Gr\"obner
cells for $\C[x,y]$ and entirely local ones for $\C[[x,y]]$.
Then we will adjust the construction 
from \cite{CV} for the former ring to the
latter. 

\subsection{Basic definitions}
As in the Introduction,
$H\!ilb^{(n)}(\C^2)$ is defined as the scheme of ideals 
$I\subset \C[x,y]$
such that $\dim \C[x,y]/I=n$. We have two standard lexicographic
orderings:
\begin{align}
&\{y^\infty\!<\!x\}:\ 1<y<y^2<\cdots <x<xy<xy^2<\cdots,\label{ord1}\\
&\{x^\infty\!<\!y\}:\ 1<x<x^2<\cdots <y<yx<yx^2<\cdots.  \label{ord2}
\end{align}
For any $f\in \C[x,y]$ let $f^0$ be its {\it maximal\,} monomial
$x^ay^b$ with respect to $\{y^\infty\!<\!x\}$. We will mainly 
need below 
$f_0$ defined as the {\em minimal\,} monomial of $f$ with respect
to $\{x^\infty\!<\!y\}$ from (\ref{ord2}). Obviously,
\begin{align}\label{min-max}
(x^ay^b+ y^{b+1}f)^0=x^ay^b=(x^ay^b+y^{b+1}f)_0 
\hbox{\ \,if\,\ } \deg_x(f)<a
\end{align}
for any  $f\in \C[x,y], a,b\ge 0$. 

For the ideals $I\subset \C[x,y]$, we set:
\begin{align}\label{defI}
I^0\equal\{f^0\mid f\in I\} \and
I_0\equal\{f_0\mid f\in I\},
\end{align}
which are {\em monomial ideals\,} by construction,
which means, as in the Introduction, that
they are linearly generated by $x^a y^b$. 

An arbitrary monomial ideal coincides with
one of $I_\la$ for the partition $\la \vdash n$
defined as follows.  Let
$\la=\{m_1\ge m_2 \ge \cdots
\ge m_\ell>0\}$, where $\sum_{i=1}^l m_i=n$. 
We set $\ell(\la)=\ell$, which is
called the length of $\la$, and $m(\la)\equal m_\ell$. 
As in the Introduction:
\begin{align}\label{Ila}
&\la'\equal
\bigl\{\{i,j\}\in \Z_+^2 \mid 0\le i <\ell,\, 0\le j < 
m_{i+1}\bigr\},\notag\\
&\hbox{and\,\,} I_\la= \oplus_{a,b}\, \C x^a y^b, 
\hbox{\,\, where\,\, }
\{a,b\}\in \Z_+^2\setminus \la'.
\end{align} 
Equivalently, 
$x^a y^b\in I_\la$ if and only if  
$a\ge i^{\circ}$ and $b\ge j^{\circ}$ for at least one {\em corner\,}
$\{i^{\circ}, j^{\circ}\}$ of $\Z_+^2\setminus \la'$. 

Any $I\in H\!ilb^{(n)}(\C^2)$ must contain pure
polynomials $f(x)$ and 
$g(y)$ in terms of $x$ and  $y$ 
of (nonzero) degree no greater than $n$. The ideals $I$ 
containing
the monomials  $x^n$ and $y^n$ form the 
{\em (local) punctual Hilbert scheme\,} $H^{(n)}$. 
More systematically:

\begin{definition}\label{defH}
The punctual Hilbert scheme is a subscheme formed by
$I\in H\!ilb^{(n)}(\C^2)$, satisfying one of the 
following equivalent conditions:

(a) $I$ contains $x^N,y^N$ for
sufficiently large $N$, which implies that $N=n$ can be
taken,  

(b) $H^{(n)}=\{I\in H\!ilb^{n)}(\C^2)\mid 
\mathfrak{m}^n\subset I\}$  for the maximal
ideal $\mathfrak{m}\equal
x\C[x,y]+y\C[x,y]$ of $(0,0)$ in $\C[x,y]$. \sq
\end{definition}

The actual {\em conductor\,} $C(I)$  of such $I$, defined
as the greatest 
monomial ideal it contains,  can be of course larger 
than $\mathfrak{m}^n$.
So only a lower bound for $C(I)$ is in this definition.
For the sake of completeness, let us check the equivalence
of $(a)$ and $(b)$.

First of all, $I\in H^{(n)}$ contains a sufficiently large
power of $\mathfrak{m}$ and can be considered naturally as
a module over $\C[[x,y]]$, which we will do constantly.
Note that the definition of $I_0$ is compatible with the
passage to $\C[[x,y]]$, since we take the smallest
monomials here.  The ring $\C[x,y]/I$ is local
artinian and the image of $\mathfrak{m}_I$ of $\mathfrak{m}$
in this ring is its maximal ideal. If $\mathfrak{m}_I^{k+1}=
\mathfrak{m}_I^{k}$ for some $k>0$, 
then $\mathfrak{m}_I^{k}=\{0\}$ by the Nakayama Lemma. 
However, such a repetition must occur no later than at  
$k=n$, since the length of the chain of consecutive 
$\mathfrak{m}_I^k$ cannot be greater than $\dim \C[x,y]/I=n$. 
\vskip 0.2cm

Similarly, let us provide the following lemma and its
justification: we want to clarify in full here the 
relation between
"global" and "local".

\begin{lemma}\label{Lemma1} 
The ideals $I^0$ for any $I\in H\!ilb^{(n)}(\C^2)$ 
and $I_0$ for any $I\in H^{(n)}$ belong to 
$H^{(n)}$, i.e. $\dim\, \C[x,y]/I^0=n$ and, correspondingly,
$\dim\, \C[x,y]/I_0=n$.
\end{lemma}  
{\it Proof.}
Let $I_\la$ be $I^0$ or $I_0$.
First of all, if $I^0=I_\la$, then any nonzero linear combination of 
$x^i y^j$ for $\{i,j\}\in \la'$
is nonzero modulo $I$; otherwise $I^0$ would
correspond to a smaller partition. Here the (linear) ordering 
of monomials can
be arbitrary, so the same argument works 
if $I_0=I_\la$ for any $I\in H\!ilb^{(n)}$, 
not only those from $H^{(n)}$. Thus $\dim\, \C[x,y]/I^0\le n$,
and the same holds for $I_0$. Let us check that such 
$\{x^i y^j\}$ linearly generate $\C[x,y]$ modulo $I$.

For the ordering $\{y^\infty\!<\!x\}$,
any monomial can be represented modulo $\{x^i y^j\}$ above
and $I$ as a sum of strictly smaller monomials. Then we continue
by induction.

In the case of $\{x^\infty\!<\!y\}$ and $I_0$,
it will be a sum of strictly bigger monomials, 
so we need an additional argument. Namely, we use
that the condition $I\in H^{(n)}$ implies that  all 
sufficiently big monomials belong to $I$. \sq

\vskip 0.2cm
The following proposition is actually a reformulation 
of the lemma.

\begin{proposition}\label{basisxy}
(i) Let $I^0=I_{\la}$ for  an arbitrary  $I\in H\!ilb^{(n)}(\C^2)$ 
or $I_0=I_{\la}$ for an arbitrary $I\in H^{(n)}$. 
Then $|\la|=n$, and the 
images of $x^i y^j$ for  
$\{i,j\}\in \la'$ form a basis of $\,\C[x,y]/I$.

(ii) For $I_\la$ as in (i), let 
$\{i^{\circ}, j^{\circ}\}$ be the corners of 
$\Z_+^2\setminus \la'$. Then $I$ is generated as an ideal
by the elements $f_{i^{\circ}\! j^{\circ}}$ such that
$f_{i^{\circ}\!j^{\circ}}-x^{i^{\circ}}y^{j^{\circ}}$
is a linear combination of $x^i y^j$ for $\{i,j\}\in \la'$.

(iii) Moreover, 
$f_{i^{\circ}\! j^{\circ}}$ for $i^\circ,j^\circ$
from $(ii)$ is unique such; it 
can contain $x^i y^j$ only if $i< i^{\circ}$ or if
$\,i=i^\circ \& j<j^{\circ}$.  
Furthermore, here  $i<i^\circ$ and $j>j^{\circ}$ must hold
in the case of $I_0=I_{\la}$ for $I\in H^{(n)}$.
\sq
\end{proposition} 

The definition of $I_0$ is compatible with the completion of
the ideals at $(x=0,y=0)$. To see this, let $\tilde{I}$ and 
$\tilde{I}_0$
be the completions with respect to $\mathfrak{m}$
of an ideal $I\subset \C[x,y]$ and $I_0$ naturally embedded 
into $\C[[x,y]]$. Here  $\tilde{I}_0$
is simply $I_0$ where infinite sums of its monomials are
allowed.  Then 
$(\,\tilde{I}\,)_0$, which is defined by picking the smallest
monomials for the same ordering
(\ref{ord2}), coincides  with $\tilde{I}_0$. 

Since  $\dim\, \C[x,y]/I\!=\!\dim\, \C[[x,y]]/\tilde{I}$
if and only if $I\in H^{(n)}$, we obtain  the following
reformulation of Definition \ref{defH}.

\begin{lemma}\label{Lemma2}
Among all ideals $I$ in $\C[[x,y]]$, the ideals $I\in H^{(n)}$ 
are characterized by the condition $I_0\in H^{(n)}$,
which is the relation $\dim C[x,y]/I_0=n$. \sq
\end{lemma}

\begin{proposition}\label{zerohn}
(i) Let us assume that an ideal $I\in \C[x,y]$ of finite 
codimension
has generators $\{f_i\}$ (as an ideal) such that 
$(f_i)^0=(f_i)_0$ and 
$I_0$ is the linear span of
$=\{(x^ay^bf_i)_0 \mid a,b\ge 0,\forall i\}$.
Then $I^0=I_0$. 

(ii) When $I\in H^{(n)}$ the generators
 $\{f_{i^\circ\!j^\circ}\}$ from Proposition \ref{basisxy}
satisfy the conditions
for $\{f_i\}$ above. Thus $I^0=I_0$ for such $I$. Vice versa, 
the conditions from $(i)$ imply that $I\in H^{(n)}$.  
\end{proposition}
{\it Proof of $(i)$}. Obviously $(x^ay^bf_i)^0=x^ay^b(f_i)^0=
(x^ay^bf_i)_0$ for any $a,b$, so $I_0\subset I^0$. The
problem can only be with linear combinations
$g=\sum c_{ab}^ix^ay^bf_i\in I$, which potentially can have $g^0$ 
smaller than $\max\{x^ay^bf_i\mid c_{ab}^i\neq 0\}$    
with respect to $\{y^\infty\!<\!x\}$.  Let
$g=x^i y^j+\sum c_{ab}x^a y^b$ with $g^0=x^iy^j$
and nonzero $c_{ab}$. 
Then either $a<i$ for any $j$, or $a=i\, \&\,  b<j$.
One has $g_0= (g')_0$ for $g'=x^i y^j+
\sum  c'_{ab}x^a y^b$, where $c'_{ab}=c_{ab}$ when
$b<j\, (\forall j)$ or $b=j\,\&\, a<i$, and $c'_{ab}$ is
zero otherwise; use
the definition of $g_0$. 
Intersecting the inequalities for $a,b$, we obtain that
$g'=\sum_{a\le i, b\le j} c_{ab}' x^a y^b$, where $c_{ij}'=1$;
see (\ref{min-max}). 
Therefore if $g^0=x^iy^j$ does not belong to $I_0$, then 
$(g')_0\not\in I_0$, which is a contradiction.
We use here that if $\{i,j\}\in \la'$ for $I_0=I_\la$,
then the whole rectangle from $\{0,0\}$ to  $\{i,j\}$ belongs
to $\la'$.
\vskip 0.2cm

{\it Proof of $(ii)$}. Using $(i)$ and Lemma \ref{Lemma2}, we 
obtain that
$I_0=I^0$ combined with $\dim\, C[x,y]/I^0=n$ gives
that $I_0\in H^{(n)}$. Without using $(i)$, the 
direct reasoning is as follows. 

Let $I_0=I_\la$ and $f=f_{0m}=y^m$ for $m=m(\la)$ in
the notations from Proposition \ref{basisxy}. We assume that
$n>0$, so $m>0$.  The partition
$\la$ is then nonempty due to part $(i)$. Indeed, it is
empty only if at least one of $f_i$
in any system of generators has a nonzero constant term. 
However this is impossible unless  $I=\C[x,y]$ due to the 
condition $(f')_0=(f')^0$. Now let us take the generator
$g$ such that $(g)_0=x^\ell$; it must exist. Here
$g=x^\ell +yp(x,y)$ and $\deg_x p<\ell$
due to $(g)_0=g^0$, but we will not need this inequality.
Then $y^{m-1}g=x^\ell y^{m-1} \!\!\mod (y^m)$ and 
$x^\ell y^{m-1}\in I$. Next, $x^{\ell}y^{m-2}g=x^{2\ell}y^{m-2}
\!\!\mod (x^\ell y^{m-1}$ gives  $x^{2\ell}y^{m-2}\in I$, and
so on. Thus $x^{m\ell}\in I$ and we can use 
Definition \ref{defH}. \sq

Finally, for partitions $\la \vdash n$, 
{\em Gr\"obner schemes} are:
\begin{align}\label{Grcells}
&Gr_\la\,=\,\{I\in H\!ilb^{(n)}(\C^2) \mid I^0=I_\la\},\ \ 
Gr_\la^0\ =\ Gr_\la\cap H^{(n)},\\
&C_\la=\{I\in H^{(n)} \mid I_0=I_\la\}\ \cong\ 
\{\tilde{I}\subset \C[[x,y]] \mid \tilde{I}_0=
(I_\la)\tilde{\,}\,\}.\label{GrC}
\end{align}

In (\ref{GrC}), we identify ideals $I\subset \C[x,y]$ with their
completions $\tilde{I}\subset\C[[x,y]]$ and $H^{n}$ with 
a scheme of ideals $\tilde{I}\subset \C[[x,y]]$ of codimension $n$.
Note that the relation $\tilde{I}_0=(I_\la)\tilde{\,}\,$, where the 
latter is the completion of $I_\la$, automatically
results in  $\dim \C[[x,y]]/\tilde{I} =n$. This 
makes the definitions of $H^{(n)}$ and $C_\la$ entirely local, 
canonically equivalent to the ones in terms of $\C[x,y]$. 

Proposition \ref{zerohn}, provides that $I_0=I^0$ for
$I\in H^{(n)}$; moreover, $\dim C[x,y]/I_0=n$ 
for $I\in H\!ilb^{(n)}(\C^2)$ results in
$I\in H^{(n)}$ due to Lemma \ref{Lemma2}. 
{\it We obtain that $Gr^0_\la=C_\la$ for any 
partition $\la$.} This is somewhat unexpected because of 
quite different definitions of $(f)_0$ and $(f)^0$.
So it  suffices to use only $C_\la$, which
we will do from now on.

\subsection{Two examples of C-cells}. Let us provide a direct 
calculation of $C_\la$ in a typical example. Generally,
the machinery of {\em syzigies\,} can be used here; 
see \cite{ES,CV,KR}. We take $\la=\{3,3,2,1\}$;
i.e. it is of order $|\la|=9$, of length 
$\ell(\la)=4$ and with $m(\la)=3$. 
The  monomials associated with the corresponding boxes
of $\la'$ are shown in Figure \ref{9-diag}. The 
monomials  without framing
are for the {\em corners\,} of $\Z^2\setminus \la'$.

\begin{figure}
{\noindent
\thicklines
%{\makebox(20,20){}}
{\makebox(20,20){}}
{\makebox(20,20){}}
{\makebox(20,20){}}
{\makebox(20,20){}}
{\makebox(20,20){j=0}}
{\makebox(20,20){j=1}}
{\makebox(20,20){j=2}}
{\makebox(20,20){j=3}}
{\makebox(20,20){}}
{\makebox(20,20){}}
{\makebox(20,20){}}
{\makebox(20,20){}}\\
\thicklines
{\makebox(20,20){}}
{\makebox(20,20){}}
{\makebox(20,20){i=0}}
{\framebox(20,20){$1$}}
{\framebox(20,20){$y$}}
{\framebox(20,20){$y^2$}}
{\makebox(20,20){$y^3$}}
{\makebox(20,20){}}
{\makebox(20,20){}}
{\makebox(20,20){}}\\
\thicklines
{\makebox(20,20){}}
{\makebox(20,20){}}
{\makebox(20,20){i=1}}
{\framebox(20,20){$x$}}
{\framebox(20,20){$xy$}}
{\framebox(20,20){$xy^2$}}
{\makebox(20,20){}}
{\makebox(20,20){}}
{\makebox(20,20){}}
{\makebox(20,20){}}\\
\thicklines
{\makebox(20,20){}}
{\makebox(20,20){i=2}}
{\framebox(20,20){$x^2$}}
{\framebox(20,20){$x^2y$}}
{\makebox(20,20){$x^3y$}}
{\makebox(20,20){}}
{\makebox(20,20){}}
{\makebox(20,20){}}\\
\thicklines
{\makebox(20,20){i=3}}
{\framebox(20,20){$x^3$}}
{\makebox(20,20){$x^3y$}}
{\makebox(20,20){}}
{\makebox(20,20){}}
{\makebox(20,20){}}\\
{\makebox(20,20){i=4}}
{\makebox(20,20){$x^4$}}
{\makebox(20,20){}}
{\makebox(20,20){}}
{\makebox(20,20){}}
{\makebox(20,20){}}
\vskip -0.4cm
\caption{$\la=\{3,3,2,1\}$}
\label{9-diag}
}
\vskip -0.3cm
\end{figure}

Accordingly, the ideals $I\in C_{\la}$, which are $I\subset H^{(9)}$
such that $I_0=I_{\la}$, are generated (as ideals in $\C[x,y]$)
by the polynomials:

\begin{align}\label{ex-9}
&f_1=x^4+C_{21}^1x^2y+C_{11}^1x^1y+C_{01}^1y
+C_{12}^1xy^2+C_{02}^1y^2,\\
&f_2=x^3y+C_{12}^2xy^2+C_{02}^2y^2,\hbox{\ and\ \, }
f_3=x^2y^2,\ f_4=y^3,\notag\\
&\hbox{where we set\ \ } f_1=f_{03},\ \,f_2=f_{22},\ \, 
f_3=f_{31},\ \, f_4=f_{40} \notag
\end{align}
in the notation $f_{i^\circ \!j^\circ}$ from
Proposition \ref{basisxy}. There are $7$ $C$-parameters
here, but $dim C_\la= |\la|-\ell(\la)=9-4=5$ due to
Theorem \ref{thmc} below. So there must be $2$ relations.
Let us find them. One has:

\begin{align*}
&yf_1= x^4y+C_{11}^1xy^2+C_{01}^1y^2 \!\!\mod (f_3,f_4),\ \
xf_1= x^5+C_{21}^1 x^3y\\
+&\,C_{11}^1x^2y+C_{01}^1xy +C_{02}^1 xy^2
\!\!\mod (f_3),\ xf_2=x^4y+C_{02}^2xy^2 \!\!\mod (f_3),\\
&\hbox{and\, \,} yf_1-xf_2=
(C_{11}^1-C_{02}^2) xy^2+C_{01}^1y^2 \!\!\mod (f_4).
\end{align*}
In the last binomial, $y^2$ is the minimal monomial. Thus
$C_{01}^1=0$, since $y^2$ belongs to (the boxes of) $\la'$. Then
$C_{11}^1-C_{02}^2$ must vanish too,
since $xy^2$ belongs to $\la'$. Finally, the relations are:
 $C_{01}^1=0=
C_{11}^1-C_{02}^2$.
\vskip 0.2cm

The following example is the most involved for the partitions
with $|\la|=5$.
Let $\la=\{4,1\}$. From  Figure \ref{5-diag}, we obtain
the following generators for any $I\in C_\la$:

\begin{figure}
{\noindent
\thicklines
%{\makebox(20,20){}}
{\makebox(20,20){}}
{\makebox(20,20){}}
{\makebox(20,20){}}
{\makebox(20,20){}}
{\makebox(20,20){j=0}}
{\makebox(20,20){j=1}}
{\makebox(20,20){j=2}}
{\makebox(20,20){j=3}}
{\makebox(20,20){j=4}}
{\makebox(20,20){}}
{\makebox(20,20){}}
{\makebox(20,20){}}
{\makebox(20,20){}}\\
\thicklines
{\makebox(20,20){}}
{\makebox(20,20){}}
{\makebox(20,20){i=0}}
{\framebox(20,20){$1$}}
{\framebox(20,20){$y$}}
{\framebox(20,20){$y^2$}}
{\framebox(20,20){$y^3$}}
{\makebox(20,20){$y^4$}}
{\makebox(20,20){}}
{\makebox(20,20){}}
{\makebox(20,20){}}\\
\thicklines
{\makebox(20,20){}}
{\makebox(20,20){}}
{\makebox(20,20){i=1}}
{\framebox(20,20){$x$}}
{\makebox(20,20){$xy$}}
{\makebox(20,20){}}
{\makebox(20,20){}}
{\makebox(20,20){}}
{\makebox(20,20){}}
{\makebox(20,20){}}
{\makebox(20,20){}}\\
\thicklines
{\makebox(20,20){}}
{\makebox(20,20){i=2}}
{\makebox(20,20){$x^2$}}
{\makebox(20,20){}}
{\makebox(20,20){}}
{\makebox(20,20){}}
{\makebox(20,20){}}
{\makebox(20,20){}}
{\makebox(20,20){}}
\vskip -0.4cm
\caption{$\la=\{4,1\}$}
\label{5-diag}
}
\vskip -0.3cm
\end{figure}

\begin{align}\label{ex-5}
&f_1=f_{20}=x^2+C_{01}^1y+C_{02}^1y^2+C_{03}^1y^3,\\
&f_2=f_{11}=xy+C_{02}^2y^2+C_{03}^2y^3,\ 
f_3=f_{04}=y^4.\notag
\end{align}

One has: $I\ni y^2f_2=xy^3 \!\!\mod (f_3)$; so $xy^3\in I$.
Using this, $xyf_2=x^2y^2 \!\!\mod (f_3, xy^3)$, and
therefore 
$x^2y^2\in I$. Now:
$$y^2 f_1=x^2y^2+C_{01}^1y^3+C_{02}^1y^4=C_{01}^1 y^3 \!\!\mod
(x^2y^2, f_3).$$
We conclude that $C_{01}^1=0$, since $y^3$ belongs to 
(the boxes of) $\la'$. Using now that $C_{01}^1=0$, we 
arrive at the $2\times 2$-system
for $xy^2,y^3$:
\begin{align*}
yf_2=\,&xy^2+C_{02}^2y^3 \!\!\mod (f_3),\\
xf_2-yf_1=C_{02}^2\,&xy^2 -C_{02}^1y^3 \!\!\mod (f_3, xy^3).
\end{align*}
It gives that $C_{02}^1+(C_{02}^2)^2=0$; otherwise 
$y^3$ would belong to $I$. 
Summarizing, $C_\la$ is obtained from
$\C^5$ by imposing  $C_{01}^1=0, 
C_{02}^1=-(C_{02}^2)^2$. So it is an affine space  of dimension 
$3=|\la|-\ell(\la)$.

It is of interest to calculate the {\em conductor\,}
$C(I)$ of $I$.  One has $x^2 f_1=x^4+C_{02}^1x^2y^2 \in I$ due
to $C_{01}^1=0$. Since $x^2y^2\in I$, we obtain that $x^4\in I$.
Finally, the generators of $C(I)$ are $y^4,xy^3,x^2y^2,x^3y,x^4$,
i.e. it is $\mathfrak{m}^4$ for generic $C$-parameters. 
When $C_{02}^1=0=C_{02}^2$, it 
will also contain $xy$ for $C_{03}^2=0$, and $x^2$ for 
$C_{03}^1=0$.
\vskip 0.2cm

The disadvantage of this direct approach is that it does not
generally provide that $C_\la$ are affine spaces. However it
can be used in any ranks and for arbitrary
isolated surface singularities, 
at least those with local rings belonging to $\C[[x,y]]$,
where Gr\"obner schemes can be readily defined. Actually,
the number of variables can be here greater than $2$.  
\vskip 0.2cm

{\sf Another example.} The following example will be needed
later; it is actually "simpler" than the one before. 
Let $\la=\{4,2\}$. From  Figure \ref{6-diag}, we obtain
the following generators of $I\in C_\la$:

\begin{figure}
{\noindent
\thicklines
%{\makebox(20,20){}}
{\makebox(20,20){}}
{\makebox(20,20){}}
{\makebox(20,20){}}
{\makebox(20,20){}}
{\makebox(20,20){j=0}}
{\makebox(20,20){j=1}}
{\makebox(20,20){j=2}}
{\makebox(20,20){j=3}}
{\makebox(20,20){j=4}}
{\makebox(20,20){}}
{\makebox(20,20){}}
{\makebox(20,20){}}
{\makebox(20,20){}}\\
\thicklines
{\makebox(20,20){}}
{\makebox(20,20){}}
{\makebox(20,20){i=0}}
{\framebox(20,20){$1$}}
{\framebox(20,20){$y$}}
{\framebox(20,20){$y^2$}}
{\framebox(20,20){$y^3$}}
{\makebox(20,20){$y^4$}}
{\makebox(20,20){}}
{\makebox(20,20){}}
{\makebox(20,20){}}\\
\thicklines
{\makebox(20,20){}}
{\makebox(20,20){}}
{\makebox(20,20){i=1}}
{\framebox(20,20){$x$}}
{\framebox(20,20){$xy$}}
{\makebox(20,20){$xy^2$}}
{\makebox(20,20){}}
{\makebox(20,20){}}
{\makebox(20,20){}}
{\makebox(20,20){}}
%{\makebox(20,20){}}
{\makebox(20,20){}}\\
\thicklines
{\makebox(20,20){}}
{\makebox(20,20){i=2}}
{\makebox(20,20){$x^2$}}
{\makebox(20,20){}}
{\makebox(20,20){}}
{\makebox(20,20){}}
{\makebox(20,20){}}
{\makebox(20,20){}}
{\makebox(20,20){}}
\vskip -0.4cm
\caption{$\la=\{4,2\}$}
\label{6-diag}
}
\vskip -0.3cm
\end{figure}

\begin{align}\label{ex-6}
&f_1=f_{20}=x^2+C_{11}^1 xy+C_{01}^1y+C_{02}^1y^2+C_{03}^1y^3,\\
&f_2=f_{11}=xy^2+C_{03}^2 y^3,\ 
f_3=f_{04}=y^4.\notag
\end{align}

One has: $I\ni y f_2=xy^3 \!\!\mod (f_3)$; so $xy^3\in I$.
Using this, $x f_2= x^2y^2 +C_{03}^2 x y^3 \!\!\mod (f_3, xy^3)$, so
$x^2y^2\in I$. Now:
$$y^2 f_1=x^2y^2+
C_{01}^1y^3+C_{11}^1 xy^3+C_{02}^1 y^4=C_{01}^1 y^3
\!\!\mod
(xy^3, x^2y^2, f_3).$$
We conclude that $C_{01}^1=0$, since $y^3$ belongs to 
(the boxes of) $\la'$, and that $C_\la$ is an affine space  
of dimension $4=|\la|-\ell(\la)$.

\subsection{The parametrization}
It is generally not true 
that $\{f_{i^{\circ}\!j^{\circ}}\}$
from  Proposition \ref{basisxy} constitute a
minimal set of generators of $I$ as an ideal. The simplest
example is \Yboxdim7pt $\yng(2,1)$. One has: $f_1=x^2+cy,
f_2=xy, f_3=y^2$, and $yf_1-xf_2=cy^2$.
Even if they are such, the number of the
corresponding coefficients of $x^i y^j$ in their decompositions
(for all corners) is generally significantly  greater than the 
dimension of the Gr\"obner cells $Gr_\la, C_\la$;
they are affine spaces. Let us address this.
  
Following the parametrization of $Gr_\la$ from \cite{CV},
we will provide an explicit parametrization of $C_{\la}$.
{\it A priori\,}, the definition of $C_\la$ is very different
from that for $Gr_\la$:\, (\ref{ord2}) is used instead of 
(\ref{ord1}) and the leading term in $f_0$ is the
{\em minimal\,} monomial, not the maximal one as in $f^0$.
Furthermore,  $C_\la$ of dimension $n-\ell(\la)$ is embedded
in $Gr_\la$ of dimension $n+m(\la)$; changing $-\ell(\la)$
to $+m(\la)$ is some combinatorial challenge too.
Theorem 1.1 from \cite{ES}
provides the dimensions
of cells in $H\!ilb^{(n)}(\C^2)$ and $H^{(n)}$, but not 
the embedding above. They use \cite{BB}, which approach
generally does not provide  explicit embeddings
of subschemes.  
Following Section 3.1 from \cite{KR}, let us
reproduce the description of  $Gr_\la$ from \cite{CV}.
\vskip 0.1cm

Given $\la=\{m_1\ge m_2\ge \cdots \ge m_\ell>0\}$ such that
$|\la|=\sum_{i=1}^\ell m_i=n$, we
set $d_1=m_\ell, d_2=m_{\ell-1}-m_{\ell},\cdots, 
d_\ell=m_1-m_2$. I.e. nonzero $d_i$ are the lengths of
the horizontal segments in the corresponding Young diagram
starting with the bottom. The $d$-set in Figure \ref{9-diag}
is $\{1,1,1,0\}$; it is $\{1,3\}$ in Figure \ref{5-diag}.
The construction is in terms of the following polynomials
in $y$: 
\begin{align}\label{pdata}
&\{p_i(y), 1\le i\le \ell \mid deg\, p_i <d_i\},\ \and \\
&\{p_{i,j}(y), 1\le i\le j\le \ell \mid deg\, p_{i,j}<d_i\}.    
\notag
\end{align}
Their coefficients will be the free parameters of $Gr_\la$. 
The following matrix of size $(\ell+1)\times \ell$ is the key:

\begin{equation*}
\t_{\la} =
\begin{pmatrix}
y^{d_1}+p_1 & 0           & 0 &\cdots & 0 & 0 \\
p_{1,1}-x   & y^{d_2}+p_2 & 0 &\cdots & 0 & 0 \\
p_{1,2}     & p_{2,2}-x   & y^{d_3}+p_3 &\cdots & 0 & 0 \\
\vdots      & \vdots      & \vdots & \ddots
 & \vdots & \vdots  \\
p_{1,\ell-2}& p_{2,\ell-2} & p_{3,\ell-2} &\cdots & 
y^{d_{\ell-1}}\!+p_{\ell-1} & 0 \\
p_{1,\ell-1}& p_{2,\ell-1} & p_{3,\ell-1} &\cdots & 
p_{\ell-1,\ell-1}\!-\!x & y^{d_{\ell}}+p_{\ell} \\
p_{1,\ell}  & p_{2,\ell}   & p_{3,\ell}   &\cdots & 
p_{\ell-1,\ell} & p_{\ell,\ell}-x 
\end{pmatrix}.
\end{equation*}

\begin{theorem}\cite{CV}.\label{thmgr}
Given $\la \vdash n$ and an arbitrary set
of polynomials from (\ref{pdata}),
the ideal $I$ generated be the  $\ell\times \ell$-minors 
of $\t_\la$ belongs to $Gr_\la$. Any $I\in Gr_\la$ can be
represented in this form for a unique set of $p$-polynomials. 
In particular,  $Gr_\la$ is an
affine space of dimension $|\la|+m(\la)=n+m_1$.\sq
\end{theorem}

The reference is \cite{CV}, Theorem 3.3.
The next theorem is the case $(i=2)$ of this theorem.
This reduction is of importance to us; we
adjust it to what we will need and
provide its complete justification (mostly following \cite{CV}).

To prevent a possible confusion while comparing
this theorem with that in  Section 3.1 of \cite{KR},
let us calculate the dimension of $Gr_\la$. It is
$$
(\ell+1)d_1+(\ell)d_2+\ldots+ 2d_\ell=\sum_{i=1}^\ell 
(\ell-i+1)d_i +\sum_{i=1}^\ell d_i= n+m_1.
$$

It suffices to take here only
the minors where lines $\{\ell+1-i^\circ\}$
are removed from $\t_\la$ for the {\em corners\,}
$\{i^\circ,j^\circ\}$ of $\Z_+^2\setminus \la'$.
Given $i^\circ$,  the corresponding minor will be 
up to proportionality a unique element $f_{i^\circ\!j^\circ}
\in I$ such that
$f=x^{i^\circ}y^{j^\circ}+\sum_{i,j} C_{ij}x^iy^j$,
where $\{i,j\}\in \la'$ and
either $i<i^\circ$ or $i=i^\circ \,\&\, j<j^\circ$. 
See Proposition \ref{basisxy}, $(iii)$.

These elements form a set of generators
of $I$. However, the coefficients $C_{ij}$ are
not arbitrary at all. They must satisfy algebraic
relations to ensure that $I^0=I_\la$, which are
"resolved" in the construction of the theorem. 
The counterpart of this theorem for $C_\la$ is as follows.

\begin{theorem}\cite{CV}.\label{thmc}
(o) Given a partition $\la \vdash n$, let $p_i=0$ in 
(\ref{pdata})
and also $p_{i,i}(y=0)=0$. Moreover, for any 
segment $[a,b]$ such that $d_{a-1}\neq 0$ and 
$\{d_a=0, d_{a+1}=0,\ldots, d_b=0\}$, 
we additionally impose the relations 
$p_{a-1,j}(y=0)=0$ for $a\le j\le b$. 

(i) Then the  $\ell\times \ell$-minors 
of $\t_\la$ generate an ideal $I$ from $C_\la$.
For a corner $\{i^\circ, j^\circ\}\!\in\! \Z_+^2\setminus \la'$, 
let $T_\la^{(\ell-i^\circ+1)}$ be the minor of  $\t_\la$ where
the line $(\ell\!-\!i^\circ\!+\!1)$ is omitted. Then 
$T_\la^{(\ell-i^\circ+1)}$ is $(-1)^{i^\circ}$ times 
$f_{i^\circ\!j^\circ}=x^{i^\circ}y^{j^\circ}+
\sum_{i,j}C_{ij}\,x^iy^j$  from
Proposition \ref{basisxy}, $(iii)$. The latter is unique
in $I$ subject to $\{i,j\}\in \la'$, and the inequalities
$i<i^\circ$,  $j>j^{\circ}$.
  
(ii) Any  ideal $I\in H^{(n)}$ can be obtained this way, 
and the corresponding set of polynomials $\{p_{i,j}\}$
subject to the conditions above is
uniquely determined by $I$. In particular, 
$C_\la$ is an affine space  of dimension $|\la|-\ell(\la)=n-\ell$,
it is naturally embedded into $Gr_\la$, and its  image is
$Gr_\la^0$ defined in (\ref{Grcells}).
\end{theorem}
{\it Proof.} Let us begin with the calculation of the
dimension of the set of $p$-polynomials. It is, indeed:
$$ \sum_{d_i\neq 0} ((l-i+1)d_i-1) - |\{d_i=0\}| =n-\ell,$$
where $|\{d_i=0\}|$ counts the additional conditions
imposed when $d_i=0$. Compare with $\dim\,Gr_\la=n+m_1$;
note that $+m_1$ is "replaced" by $-\ell$.  We will set
$d_0=0$ and $p_0=0$ later on. See \cite{CV}, Corollary 3.1
$(i=2)$.

Due to Propositions \ref{basisxy} and  \ref{zerohn},
the following
property is necessary for the minors $T_\la^{(i)}$, which are
determinants of  $\t_\la$ in Theorem \ref{thmgr}
without line $i\ (1\le i\le \ell+1)$ and under the assumption
that $\ell-i+1$ is  $i^\circ$ from some corner
$\{i^\circ,j^\circ\}$ of $\Z_+^2\!\setminus\! \la'$.

The property which must hold is that
\begin{align}\label{mform}
&T_\la^{(i)}=\pm\, x^{\ell-i+1}y^{d_0+\ldots+d_{i-1}}+
\sum_{a,b}C_{ab}\,x^a y^b,\, \where\\
& \{a,b\}\in \la'\and  a<\ell\!-\!i\!+\!1, \ 
b>d_0+\ldots+d_{i-1},\notag
\end{align}
for some $C_{ab}$ if and only if the conditions 
$p_i=0$ and the other ones from $(o)$ are imposed.
I.e. that these conditions for $p_i$
 are necessary and sufficient for the inequalities for $\,a,b\,$
in (\ref{mform}).  

Recall that $\ell-i+1$ is some $i^\circ$ only
when $d_i\neq 0$ or when $i=\ell+1$. Provided $(o)$,
if $\ell-i'+1$ is not assumed 
to be $i^\circ$, then $T_\la^{(i')}$
is $\pm\, x^c$ multiplied by $T_\la^{(i)}$ 
from (\ref{mform}) for a corner $i^\circ=\ell-i+1$ such that 
$d_i\neq 0,\, d_{i+1}=0,\,\cdots,\, d_{i+c}=0$. So $c$ is
the distance from $i'$ to the greatest possible $i$ 
such that  $d_i\neq 0, i\le i'$.

Generally, 
$T^{(i)}_\la=\prod_{j=0}^{i-1}(y^{d_j}+p_j)\det \t^{[i]}_\la$,
where $\,1\le i\le \ell+1\,$,\, for

\begin{equation*}
\t^{[i]}_\la = 
\begin{pmatrix}
p_{i,i}-x   & y^{d_{i+1}}+p_{i+1} & 0 & \hspace{-3ex}
\cdots 0 & 0 \\
p_{i,i+1}     & p_{i+1,i+1}-x   &y^{d_{i+2}}+p_{i+2} &  
\hspace{-3.5ex}\cdots\,\vdots &  0 \\
\vdots     & \vdots      & \hspace{3ex}\vdots\ddots
& \hspace{-3ex}\cdots 0 & \vdots  \\
p_{i,\ell-2}& p_{i+1,\ell-2} &\cdots& 
y^{d_{\ell-1}}\!+p_{\ell-1} & 0 \\
p_{i,\ell-1}& p_{i+1,\ell-1} &\cdots & 
p_{\ell-1,\ell-1}\!-\!x & y^{d_{\ell}}+p_{\ell} \\
p_{i,\ell}  & p_{i+1,\ell}   &\cdots & 
p_{\ell-1,\ell} & p_{\ell,\ell}-x 
\end{pmatrix}.
\end{equation*}

Since the structure of these matrices is uniform with
respect to $i$, only  $\t^{[1]}_\la$ are sufficient to consider.
Let us make this exact. We set $\de_0=0$, 
$\de_i=|\{1\le j<\ell\!-\!i\!+\!1 \mid m_j=m_{\ell-i+1}\}|\,$
 for $1\le i\le \ell$. Let  $\la[i]$ be the partition given  
by $\{m_1\!-\!m_{\ell-i+1},\, 
m_2\!-\!m_{\ell-i+1},\ldots, m_{\ell-i}\!-\!m_{\ell-i+1}\}$ 
where we omit
the last $\de_i$ zeros. Geometrically, we remove the first $i$ 
columns from the diagram describing $\la$. Then
for $\,1\le i\le \ell+1$:

\begin{align}\label{indM}
&\det \t^{[i]}_\la=(-x)^{\de_{i-1}} 
\det \t^{[1]}_{\la[i-1]}, \and \\
&T^{(i)}_\la=(-x)^{\de_{i-1}}\prod_{j=0}^{i-1}(y^{d_j}+p_j)
 \det \t^{[1]}_{\la[i-1]},
\ \notag
\end{align}
where the indices of $p$\~polynomials in  $\t^{[1]}_{\la[i-1]}$
must be as in $\t^{[i]}_\la$ where the first $\de_i$ 
columns and rows are deleted;  here $\la[0]\equal \la$.

We can now proceed by induction with respect to $n=|\la|$,
where the case $n=1$ is obvious. This is helpful but not
actually needed below. 
Let us check that the conditions
from $(o)$ are necessary for $I\in H^{(0)}$.
\vskip 0.1cm

First of all, without any induction, $T_\la^{(\ell+1)}=
\prod_{j=0}^{\ell}(y^{d_j}+p_j)$ is $y^{m_1}$ if and only if
all $p_j$ are zero. Let as now assume that the entries of
$\t^{[1]}_{\la[1]}$ with the shift of indices as in (\ref{indM})
satisfy the conditions from $(o)$.
Then $p_{1,1}$ must be divisible by $y$, since otherwise
$T^{[1]}_{\la}$ will contain an $x$-monomial of degree
smaller than $\ell$ because of the contribution of the
diagonal. Here the usage of the induction is not necessary too.
\vskip 0.2cm

The last check concerns the additional conditions
addressing  the columns where $d_i=0$ for some $i$.  
Let $d_{i+1}=0$ and 
$d_{i}\neq 0$, i.e. $\ell-i+1$ is $i^\circ$ for some
corner. Then $T_\la^{(i)}$ must have the smallest monomial
$x^{\ell-i+1}y^{d_0+\ldots+d_{i-1}}$ due to (\ref{mform}).
However, if $p_{i,i+1}$ has a nonzero constant term, then it
must contain $x^{\ell-i-1}y^{d_0+\ldots+d_{i-1}}$, which
is smaller  with respect to the ordering $\{x^\infty\!<\!y\}$.
This is impossible. 

More generally, let $d_{i}\neq 0$,
$d_{i+1}=0=\cdots=d_{i+r}$, $d_{i+r+1}\neq 0$.
We set $\a^{[i]}_\la=\t^{[i]}_\la(y\!=\!0)$,
$a_{i,j}=p_{i,j}(y\!=\!0)$; we set $r=3$
as in the picture.  
Using what we have already checked, the matrix $\a^{[i]}_\la$ is as
follows:

\begin{equation*}
\hspace {-20ex}\a^{[i]}_\la = 
\begin{pmatrix}
-x            & 1    & 0  & 0  & 0 &\cdots & 0 \\
a_{i,i+1}     & -x   & 1  & 0  & 0 &\cdots & 0 \\
a_{i,i+2}     & 0    & -x & 1  & 0 &\cdots & 0 \\
a_{i,i+3}     & 0    & 0  & -x & {\bf 0} &\cdots & 0 \\
a_{i,i+4}     & 0    & 0  & 0  & -x &\cdots & 0 \\
\vdots        & \vdots &\vdots       &\vdots & 
\vdots & \ddots & \vdots \\
a_{i,\ell}    & 0 & 0 & 0  
&a_{i+4,\ell}   & \cdots & -x 
\end{pmatrix}.\hspace{-41ex}{\linethickness{0.1pt}
\framebox(110,55){}}
\end{equation*}

The (bold) entry $\{4,5\}$ equals $0$ because it
is $y^{d_{i+r+1}}\!=\!y^{d_{i+4}}$ at $y\!=\!0$. Thus, the 
determinant of $\a^{[i]}_\la$ is the product of
the determinant of its upper principal $4\times 4$\~block and the
determinant of the principal block starting with the entry 
at $\{5,5\}$.
The former determinant is classical: 
$x^4-(a_1x^2+a_2x+a_3)$ for $a_j=a_{i,i+j}$. 
Unless $a_1=a_2=a_3=0$, the final product cannot be
a pure $x$\~monomial, as it is supposed to be.
This gives the required.  
 
These arguments can be equally  used in the opposite 
direction; they actually give that conditions $(o)$ are not only
necessary but sufficient too, i.e. equivalent to the 
summation restrictions in (\ref{mform}). A direct deduction of 
these restrictions from $(o)$ is not difficult as well. 
Let us emphasize that
we rely  in this
proof on Theorem \ref{thmgr}, which provides that no polynomials 
strictly within (the boxes of)
$\la'$ can occur in $I$ generated by $T_\la^{(i)}$. 
Since it holds for $\t_\la$, then of course it
is true under any specializations of the coefficients
of $p$\~polynomials. 
\sq

\section{\sc Plane curve singularities}
\subsection{Compactified Jacobians}\label{sec:CompJ}
We will define their Gr\"obner decomposition. Let us
consider only unibranch plane curve singularities. They
are  subrings $\r\subset \C[[z]]$ in terms of
the uniformizing parameter $z$ that have 2
generators, $x$ and $y$, and 
have $\C((z))$ as the field of fractions. We set
$\nu_z(uz^a+\sum_{i>a}c_iz^i)=a$ for $u\neq 0$,
 and define the {\em valuation
semigroup\,} $\Ga=\Ga_\r=\{\nu_z(f)\mid f\in \r\}$. 
See here and below
\cite{PS,GP} or \cite{Ch,ChP1}. One has:
 \begin{align}\label{rval}
\de\equal \dim_{\C} \C[[z]]/\r = |\Z_+\setminus \Ga|\,, \and 
\r\supset (z^{2\de}).
\end{align}
From now on, $(z^m)\equal\C[[z]]z^m$ for $m\in \Z_+$;
it is a principal ideal in $\C[[z]]$, not in $\r$.

The {\em compactified Jacobian\,} of $\r$, 
denoted by  $\overline{Jac}_\r$ or simply by $\overline{Jac}$ 
 is a projective 
reduced irreducible scheme defined as follows:
\begin{align}\label{cjac}
&\overline{Jac}\equal \{\r-\hbox{submodules}\,\, 
M\subset \C[[z]]\,\mid\,  \dim_{\C} C[[z]]/M=\de\}.
\end{align}
As in \cite{PS}, the structure of projective variety 
is due to the following:
\begin{align} \label{cjacde}
&\hbox{for any\ \,} M\in \overline{Jac}, \ 
M\supset (z^{2\de})=\C[[z]] z^{2\de}. 
\end{align}

An important invariant of an $\r$\~module $M\subset \C((z))$ is
$\De(M)\equal\{\nu_z(v)\mid v\in M\}$.
For $M\subset \C[[z]]$, which we will always assume
below unless stated otherwise,
and for $k\in \Z_+$:
$$\dim \C[[z]]/z^k M=k+\dim \C[[z]]/M,\ \dim \C[[z]]/M=
|Z_+\setminus \De(M)|. $$

Let $M_\bullet\equal 
z^{-v}M$ for $v=v(\De)\equal
\min \De$, and $\De_\bullet(M)=\De(M)-v$. The latter is
a $\Ga$\~module,
which means by definition that $\Ga+\De\subset \De$.  
One has: $\De_\bullet(M)=\De(M_\bullet)$. 

The image of the map $M\mapsto M_{\bullet}$ is the
set of all {\em standard} modules in $\C[[z]]$, which
are those containing some elements in $1+z\C[[z]]$.
So $v=0$ for such modules and $M_\bullet=M$ if and only
if $M$ is standard.
Equivalently, $M$ is standard if $M\C[[z]]=\C[[z]]$ or
if $\De[M]$ is standard, where a $\Ga$-module $\De$ is
called standard if $0\in \De\subset \Z_+$.
\vskip 0.2cm

The {\em generalized Jacobian\,} $Jac\subset \overline{Jac}$ 
is formed by
all {\em invertible
modules\,}, i.e. $M_\phi=\r \phi$ for
$\phi\in 1+z\C[[z]]$.
Equivalently, invertible modules are such that 
$\De(M)=\Ga$. 
The generalized Jacobian is a group  with respect to the  multiplication
of the generators $\phi$. It is an affine 
space of dimension $\de$; the simplest parametrization
is as follows:  
$\phi=1+\sum_g \phi_g z^g$, where $g\in \Z_+\setminus \Ga$,
$\phi_g\in \C$.
Its closure is the whole $\overline{Jac}$. Let
$Jac^{\bullet}\equal \{ M=M_\bullet\}$, i.e. 
it is formed set-theoretically by all
standard $\r$\~modules. It is a disjoint union of
quasi-projective schemes $Jac^{(d)}\equal\{M=M_\bullet \mid
\dim C[[z]]/M =d\}$; note that  $Jac=Jac^{(\de)}$. 
\vskip 0.2cm

We will also use the {\em duality\,} (also called
{\em reciprocity\,}).
 For an $\r$-module
$M\in \C((z))$, let $M^\ast\equal
\{f\in \C((z))\mid fM\in \r\}$. It is an $\r$\~module.
One has: $(M^{\ast})^\ast=M$, 
$\C[[z]]^\ast=(z^{2\de})$, and $(\r)^\ast=\r$.

Let us assume that $M=M_\bullet$,
i.e. that $M$ is {\em standard}. Then 
$(M)^\ast\in \C[[z]]$, 
since $M$ 
contains $1+z(\cdot)$, and it belongs to $\r$ if and
only if $M$ contains $\r$.
The latter gives that $\De(M^\ast)\subset \Ga$, which 
is another defining property of standard $M$:
$$
M=M_\bullet \Longleftrightarrow \De(M^\ast)\subset \Ga.
$$
Indeed, if $\De(M^\ast)\subset \Ga$ then $\Ga\subset \De(M)$,
and therefore $M$ is standard. 
\vskip 0.2cm

The following
holds for standard $M$ and for any Gorenstein rings 
$\r\subset \C[[z]]$,
not only for the rings of plane curve singularities
we consider in this paper; see \cite{PS} and \cite{GM2}.
Let $\De=\De(M)$. Then
$\De(M^\ast)=\De^\ast\equal \{p\in \Z\mid p+\De\subset \Ga\}\subset
\Z_+$
and 
\begin{align}\label{Dedual}
\De^\ast= (2\de-1)-(\Z\setminus \De)=
(2\de-1-(\Z_+\setminus \De))\cup (2\de+\Z_+),
\end{align}
where $v+X=\{v+x \mid x\in X\}$ for $X\subset \Z \ni v$.
It gives that $M^\vee\equal z^{c(M)-2\de}M^\ast$ is again a 
standard module 
for the {\em conductor\,} $c(M)=c(\De)$, which is defined as
$\min\{c\mid (z^c)\subset M\}=\min\{c\mid c+\Z_+\subset \De\}$.
We call $M^\vee$ the {\em standard dual\,} of $M$. 
The conductor is obviously  $g_{t\!op}(\De)+1$ for the top element
$g_{t\!op}$ in $\Z_+\setminus \De$.
Thus $\De((M^\ast)_\bullet)= \De^\vee\equal
(c(\De)-1 -(\Z_+\setminus \De))\cup (2\de+\Z_+)$  for a
standard $M$. We call $\De$ {\em selfdual} or, more exactly,
{\em standard selfdual}
if $\De=\De^\vee$. 
\vskip 0.1cm

{\sf Applications.}
The minimal embeddings $f M\subset \r$ are
as follows. Let
\begin{align}\label{u-min}
u_{min}=\min\{u \mid f_u M\subset \r \for f_u-z^u\in
(z^{u+1})\}.
\end{align}
Then  $u_{min}=\min \De^\ast=
2\de-1-g_{t\!op}(\De)=2\de-c(M)$
due to (\ref{Dedual}), which holds for any $\r$-modules 
$M$,
not only standard. 
\vskip 0.1cm

As another application, let $J^1\equal \{M^1\mid M^1\supset \r\}$, so 
the modules $M^1\subset \C[[z]]$ here are automatically standard;
$J^1$ is a disjoint union of
{\em projective\,} schemes with  fixed $\dim C[[z]]/M^1$.
One has: 
\begin{align}\label{UJac}
Jac^\bullet=(U^1/U_\r^1) \,J^1 \for U^1\equal 1+z\C[[z]],\, 
U_\r^1=U^1\cap \r.
\end{align}
Indeed, one can obtain
any standard $M$ as $\phi M^1$ for proper $M^1$ and 
invertible $\phi\in \C[[z]]$. So $J^1$ is a certain {\em skeleton\,}
of $Jac^\bullet$:

\begin{lemma}\label{lem:tilde-Jac}
Let \, $\tilde{Jac^\bullet}
\!\equal\!\bigl\{ (M,\phi) \mid M\!=\!M_\bullet,\,
\phi\!\in\! (U^1\cap M)/U_{\r}^1\bigr\}$. 
Then $\tilde{Jac^\bullet}$ is
isomorphic as a scheme to $J^1\times (U^1/U_\r^1)$, and 
$\tilde{Jac}^{(d)}\equal 
\{(M,\phi) \mid M\!\in\! Jac^{(d)}\}
\simeq
Jac^{(d)}\times \C^{\de -d}$ (as schemes).
\comment{ Let $J[\De]=\{M=M_\bullet \mid \De(M)=\De\}, 
J^1[\De]=J[\De]\cap J^1$, $\tilde{Jac}[\De]=
\{(M,\phi)\in \tilde{Jac^\bullet} \mid M\in J[\De]\}$.
Then $\tilde{Jac}[\De]\simeq J^1[\De]\times(U^1/U_\r^1)\simeq$
$J[\De]\times \C^{\de -d}$ for $d=|\Z_+\setminus \De|$ 
(as schemes).}
\end{lemma}
{\it Proof.} The action $(M,\phi)\mapsto f(M,\phi)=(fM,f\phi)$ 
of $f\in U^1/U_{\r}^1$ on the pairs $(M,\phi)$ is free
because it is free at the $\phi$-component. Then 
$(M,\phi)=\phi(M/\phi,1)$, where $M/\phi$ contains $1$ and
belongs to $J^1$.

To finish the proof, we pick $f_g\in \r$ for $g\in \Ga$
such that $f_g=z^g \mod (z^{g+1})$ in $\C[[z]]$.
Let $\De=\De(M)$ for a standard $M$, 
$|\Z_+\setminus \De|=d$ and 
$\De\setminus \Ga=\{g_1<g_2<\ldots<g_{\de-d}\}.$
Then any $\phi=1+\sum_g \phi_g z^g$ in $(M,\phi)$
can be represented (modulo $U_{\r}^1$) 
%modulo $(z^{2\de})\subset M$ 
as $\phi^\circ=1+\sum_{i=1}^{\de -d} 
\phi_i z^{g_i}+\sum_{g\not\in \De} \phi'_g z^g$.
We proceed here by induction. 
Let $\phi_g\neq 0$ for the minimal such $g\in \Ga\setminus\{0\}$. 
Then $\psi=\phi(1-\phi_g f_g)$ has $\psi_h=0$ for any $h\le g$ 
in $\Ga\setminus \{0\}$.

The coefficients $\phi_{i}$ determine $\phi^\circ$
uniquely. 
Indeed, $\psi_g=c$ for $\psi=\phi^\circ (1+cf_g)$, 
any $c\in \C^*$ and $g\in \Ga\setminus \{0\}$.
Any $\phi_i\in \C$ can occur here:\ $M$ modulo  
$\{\sum_{g\not\in \De} \phi'_g z^g\}$
has a basis $z^g$ for $g\in \De$.  
Thus, $\tilde{Jac}^{(d)}\simeq
Jac^{(d)}\times \C^{\de -d}$,
where $\phi_i$ are
the coordinates of the latter factor.
%$\C^{\de -d}=\{\phi_g,\ g\in \De\!\setminus\!\Ga\}$
\sq 
\vskip 0.1cm

%We will use this lemma upon the
%$\ast$-identification of $J^1$ with 
%$J_1\equal
%\{M\mid (z^{2\de})\subset M\subset \r\}$, 
%a subset of the set of all ideals in $\r$.
%\vskip 0.1cm 

{\sf Choosing $x,y$.}
Changing the parameter $z$, we can assume that the generators
of $\r$ are $y=z^a, x=z^b h(z)$,
 $h=1+\sum^c_{j=1} u_j z^j$ for some $u_j$, and $1<a<b$.
Here the total $gcd$ of $a,b$ combined with
the degrees of monomials in $h$
must be $1$. Also, $b$ can be assumed the smallest 
in $\Ga$ non-proportional to $a$, and 
$z^j$ with  $j+b\not\in \Ga$ are sufficient.
The sharp version of the latter restriction is that
$h\!=\!1\!+\!z^{\la-b}+\sum_{j>\la-b}^c u_j z^j$, where
$j+b\not\in\La'$ for 
$\la\equal\min \{\La'\setminus \Ga'\}$, 
the {\it Zariski invariant\,}, and 
\begin{align}\label{Laprime}
&\ \Ga'\equal(\Ga\setminus 0)-a\,\subset\, 
\La'\equal \nu_z\bigl(\r z^{a-1}+\r (d(z^bh)/dz)
\bigr)\!+\!1\!-\!a. 
\end{align}
The latter is the 
$\Ga$-module of "K\"ahler differentials", shifted to make it
a {\it standard} one.  See e.g. \cite{HH}. 
\vskip 0.cm

Let  $P(x,y)=res (z^a-y, z^b h(z)-x; z)$, the 
{\em resultant\,} with respect to $z$. Then $P(z^bh(z),z^a)=0$ and
$(-1)^a P=x^a +(-1)^{(a+1)(b+c)} y^{b+c}+
\sum_{i<a,j<b+c} d_{i,j} x^i y^j$, where
$i+j\ge a$, which is direct from the definition of the resultant.
Also, $(-1)^aP=
\prod_{i=0}^{a-1}\bigl(\,x-(\om^i z)^b\,h(\om^i z)\,\bigr)$ for
$\om=\exp(\frac{2\pi \imath}{a})$.
The polynomial $P$ is irreducible.

\begin{lemma}\label{Peq}
Let $P$ be as above, $d_{i,j}\neq 0$ and $i=a-v$. Then 
$$j\in \j_v\equal \bigl\{\frac{vb+r_1+\ldots+r_v}{a}\in \N 
\,\mid\, 0\le r_k\le c, \,1\!\le\! k \!\le\! v\bigr\},$$
\noindent
where $r$ are $0$ or from the set of indices such that $u_r\neq 0$
in $h(z)$. 
Setting $\ep_v\!=\!\min \j_v\!-\!vb/a\ge 0$, one has: 
$ib+ja\!=\!\nu_z(x^iy^j)\!\ge\! a(b\!+\!\ep_v)$.
\end{lemma}
{\it Proof.} Consider $x$ as an element of $\C(y^{1/a})$:\,
$x=y^{b/a}h(z\mapsto y^{1/a})$. It generates this field over
$\C(y)$ by our assumptions.
The coefficient
$d_{a-v,j}$ is then the sum of the {\em traces\,} in the
field extension $\C(y^{1/a})/\C(y)$ of the elements 
$\C$-proportional
to $y^{vb/a+(r_1+\ldots+r_v)/a}$, where $0\le r_k\le c$.
The exponent here must be an integer, which provides that
$j\ge \ep_v+vb/a$. So $ib+ja
\ge (a-v)b+(\ep_v p+vb/a)a=a(b+\ep_v)$.
\sq 
 
For instance, let us consider the  simplest unibranch
plane curve singularity that is not quasi-homogeneous,
which is for the ring
$\r=\C[[y=z^4,x=z^6+z^7]]$. Here one has: $a=4,b=6,c=1$ and
$P=x^4-2x^2y^3-4xy^5+y^6-y^7$. Indeed, $\j_1=\emptyset$, i.e.
$x^3$ does not appear, 
$\j_2=\{(12+\{0,1,2\})/4\}\,\cap\, \N=\{3\}$, which corresponds
to $x^2y^3$, 
$\j_3=\{(18+\{0,1,2,3\})/4\}\,\cap\, \N=\{5\}$, and
$\j_4=\{(24+\{0,1,2,3,4\})/4\}\,\cap\, \N=\{6,7\}$, which
gives $y^6,y^7$. The exact
calculation of the coefficients of $P$ using the
$\j$-sets is straightforward. For generic 
$\{d_{i,j}\}$, the estimates for $j$ from the Lemma can be reached.  
The structure of $P(x,y)$ is of importance in
what will follow.

\comment{
Let us involve  the ring of singularity
$\w\equal \C[[x,y]]/(\partial P/\partial x,
\partial P/\partial y)$. Then $\dim_\C \w=2\de$.
This dimension can only increase if we degenerate $P$
to  $P_0(x,y)=x^a+y^b$,
where it is equal to $(a-1)(b-1)$. 
Combining this with Lemma \ref{Peq}, we obtain that
\begin{align}\label{deabineq}
& 2\de\le (a-1)(b-1),  ib+ja\ge ab> 2\de \for C_{i,j}\neq 0,
\end{align}
i.e. all monomials $x^i y^j$ in $P$ satisfy 
$\nu_z(x^i y^j)\ge 2\de$. Actually, we arrive at a somewhat
stronger fact: $(i-1)a+(j-1)b>2\de$. Geometrically, the
box corresponding to $\{i,j\}$, which is with its
upper-left corner $\{i-1,j-1\}$ is strictly below the
anti-diagonal $ua+vb=2\de$, where the coordinates $u,v$
are positive in the lower-right quarter.    
}

Note that the monomial ideal in $\C[[x,y]]$ linearly
generated by $x^iy^j$ such that $\nu_z(x^iy^j)\ge ab$,
which contains $P(x,y)$,
maps into $(z^{2\de})$
if $2\de\le ab$; the latter
holds when $gcd(a,b)=1$ and in some other cases.

\subsection{Topological invariance}
In the classical geometry of smooth projective
curves, they can be "recovered" from 
$\overline{Jac}$ supplied with the polarization divisor.
%The submanifolds of effective divisors of degree $1$ is
%the curve itself. 
We will check that
$\tilde{Jac^\bullet}$ from Lemma \ref{lem:tilde-Jac}
is a {\em topological\,} invariant of $\r$
for some basic families.
Also, some additional structures of
$\overline{Jac}$ allow the  "extraction"
of the equation of the initial singularity from it; 
see Proposition \ref{PfromI}, (i).  

The topological type of $\r$ is given by the isotopy
type of its {\em link\,}: the intersection
of $\{(x,y)\mid P(x,y)=0\}$ with sufficiently
small $\S^3$ centered at $(0,0)$. The semigroup $\Ga$ fully
determines it, which is a classical 
fact; see \cite{Za}. The topological
invariance of $\tilde{Jac^\bullet}$ as a scheme is generally subtle.
If it holds, then it  readily gives the  topological invariance of  
the {\em motivic superpolynomial\,} of $\r$
conjectured in \cite{ChP1}; the demonstration is below.

\begin{proposition}\label{MIast}
Let $\tilde{C}_\r$ be the full preimage of $(z^{2d})\subset \r$
in $\C[[x,y]]$, i.e.  $\tilde{C}_\r=\C[[x,y]]]P(x,y)+
\C[[x,y]](z^{2\de})'$ for any set-theoretical lift of $(z^{2\de})'$
of $(z^{2\de})$ to $\C[[x,y]]$.  
We assume that for some family of rings $\r$,
$\Ga$ is fixed and 
$\tilde{C}_\r$ is constant considered up to automorphisms 
of the ring $\C[[x,y]]$.
Then the schemes 
$\tilde{Jac^\bullet}$ and $\tilde{Jac}^{(d)}$
from Lemma \ref{lem:tilde-Jac}
considered up to isomorphisms 
are constant. Here  $\tilde{Jac}^{(d)}$ 
covers $Jac^{(d)}$ for $0\le d\le g$ with the
fibers isomorphic to $\C^{\de-d}$ 
due to this lemma.
\end{proposition}
{\it Proof.} One has:
 $\tilde{Jac^\bullet}=\{(M,\phi)\mid
\phi\in (U^1/U_\r^1)\cap (M/U_\r^1)\}$ for standard $M$; it is  
a disjoint union of $\tilde{Jac}^{(d)}$ for $d=|Z_+\setminus \De(M)|$.
Here
 $U^1= 1+zC[[z]]$,
$U_\r^1=U^1\cap \r$ and  $U^1/U_\r^1\simeq \C^\de$ (as spaces). 
One has: $\tilde{Jac^\bullet}\simeq 
J^1\times (U^1/U_\r^1)\simeq 
J_1\times (U^1/U_\r^1)$, where $J^1=\bigl\{M=M_\bullet\mid 1\in M
\bigr\}$
and 
$J_1\equal\bigl\{
z^{2\de}\C[[z]] \subset M^* \subset \r\,\bigr\}$; we use the
duality map $M\mapsto M^*$, which identifies $J^1$ with  
$J_1$. By sending ideals $M^*$ to their
full lifts $\tilde{M}^*$ in $\C[[x,y]]$, we identify $J_1$ with 
the scheme $\bigl\{\text{\,ideals\,} \tilde{I} \subset \C[[x,y]] \mid
\tilde{C}_\r\subset \tilde{I}\,\bigr\}$. Such $\tilde{I}$ are
exactly $\tilde{M}^*$ by construction.
Thus, $\Ga$ and $\tilde{C}_\r$ completely determine
$\tilde{Jac^\bullet}$ (up to isomorphisms). 

Now fix
$d=|Z_+\setminus \De|$ for standard $\De=\De(M)$. 
From (\ref{Dedual}):
$\De^\ast=
(2\de-1-(\Z_+\setminus \De))\cup (2\de+\Z_+)$, and the corresponding
$d^\ast$ equals $2\de-d$. Then $\tilde{J}^1[\De]=
J^1[\De]\!\times\! (U^1/U_\r^1)\simeq 
J_1[\De^\ast]\!\times\! (U^1/U_\r^1)$ under the map $\ast$.
Here $J^1[\De]=J[\De]\cap J^1$ and $J_1[\De^*]=J[\De^*]\cap J_1$.
%recall that $(M^\ast)^\ast=M$. 
%Recall that 
%$(M^\ast)^\ast=M$, $(\De^\ast)^\ast=\De$.
Therefore
$\tilde{Jac}^{(d)}$ is  the product of $(U^1/U_\r^1)$ and
the scheme
$\bigl\{M^*\in J_1\mid |\Z_+\setminus \De(M)|=2\de-d\bigr\}=$
$\Bigl\{\tilde{C}_\r\subset \text{\,ideals\,} \tilde{I} \subset 
\C[[x,y]]\, \mid\,
 \dim \C[[x,y]]/\tilde{I}=2\de-d\,\Bigr\}$, i.e. constant. 
%It covers $Jac^{(d)}$ with the fibers isomorphic to $\C^{\de-d}$. 
We use that 
$\dim \C[[x,y]]/\tilde{M}^*=\dim \r/M^*=2\de-d$, where $\tilde{M}^*$ is
the full lift (preimage) of the ideal $M^*\subset \r$. 
\sq
\vfil

This proposition can be extended to the {\it flags\,}
of standard modules and ideals from \cite{ChP1}; we will omit
the details. 
%Actually, its justification is of general
%kind. The pair $\C[[x,y]]\subset \C[[x,z]]$ is connected 
%with the instantons of rank $a$; similar pairs occur for
%surface singularities. 
Examples of "constant"
$\tilde{C}_\r$ will be provided below. When  $\tilde{C}_\r$ is a monomial ideal in $\C[[x,y]]$,
\cite{BB} can be used to prove that
$\tilde{Jac}[\De]$ 
are affine cells.
\comment{
Later, we will discuss
an approach based on 
"stable" $(z^{2\de+e})$, which is for sufficiently large $e$. The
images of the natural embeddings 
of $\overline{Jac}, Jac^\bullet$ into
$\cup_{n} H^{(n)}$ can be described by explicit equations
for such $e$. }
\vskip 0.2cm
\vfil

{\sf Motivic superpolynomials.}
They are defined for the rings
$\r\subset \F_q[[z]]$ over $\F_q$ (the field with $q$ elements):
$\h_{\r}(q,t)\equal
\sum_{M}t^{dim_{\F_q}(\F_q[[z]]/M)}$, where the
summation is over standard $M$,
$\r$-submodules $M\subset \F_q[[z]]$ such that
$M\F_q[[z]]=\F_q[[z]]$. The ring $\r$ is supposed to be with
$2$ generators and with the field of fractions $\F_q[[z]]$.
For any $\r$ over $\C$, we can define it over $\Z$ and then
over  $\F_q$ for any $q=p^m$ for sufficiently general prime $p$
(apart from finitely many of them) with the same valuation 
semigroup $\Ga$.   The preservation of
$\Ga$ is the 
weakest possible definition of places $p$ of 
{\em good reduction}. The
$\Ga$-module of K\"ahler differentials (see above)
must be also assumed to remain unchanged upon the passage 
from $\C$ to $\F_q$ for some considerations.   
%\vfil

For $\tilde{Jac^\bullet}$, we naturally set:
$\tilde{\h}_{\r}(q,t)\equal
\sum_{(M,\phi)}t^{dim_{\F_q}(\F_q[[z]]/M)}$, where the summation
is over $(M,\phi)$ in $\tilde{Jac^\bullet}$ defined over $\F_q$.
Due to Lemma \ref{lem:tilde-Jac}, which holds for rings $\r$
over $\F_q$,
the number of $\phi$ for a given $M$ equals $q^{\de-d}$ for
$d=dim_{\F_q}(\F_q[[z]]/M)$. Thus, 
 $\tilde{\h}_\r(q,t)= q^\de \h_\r(q,t/q)$, and we see that the usage of
$\tilde{Jac^\bullet}$ is sufficient to obtain $\h_\r(q,t)$.
The number of $\F_q$-points of $\tilde{Jac}^{(d)}$ depends only
on $\Ga$ if $\tilde{C}_\r$ depends only on $\Ga$ due to 
Proposition \ref{MIast}; the latter holds over $\F_q$ (if $p$ is
a place of good reduction). Thus 
 $\h_\r(q,t)$
is a topological invariant for such $\r$. 

%Actually, here any Gorenstein rings
%$\r$ can be considered here for a proper replacement
%of $\C[[x,y]]$, but we (conjecturally) 
%need plane curve singularities
%for the functional equation for $\h_\r(q,t)$, which is under
%$q\leftrightarrow t^{-1}$.

\comment{
Let us extend $\C[[x,y]]$ by $z$ satisfying $z^a=y$;
thus $\C[[x,y]\subset \C[[x,z]]$ for $y\mapsto
z^a$. The natural projection of $\C[[x,y]]$ onto 
$\r$ can be extended to the map  
from $\C[[x,z]]$ to the product 
$\C[[z]]^a=\prod_{i=0}^{a-1}\, \{\C[[z]]\}_i$, denoted also by 
$(\{\C[[z]]\}_i)$,  of the $\,a\,$ copies of the ring $\C[[z]]$:
$$
z\mapsto (z,\ldots,z),\ y\mapsto (z^a,\ldots,z^a),\ 
x\mapsto \bigl((\om^i z)^b\,h(\om^i z)\in \{\C[[z]]\}_i\bigr)
$$
for $h=1+\sum^c_{j=1} u_jz^j$ as above and
$\om=\exp\frac {2\pi \imath}{a}$. 

The ideal 
$(x,z^b)\equal\C[[x,z]]x+\C[[x,z]]z^b$ contains $P(x,y)$
since it contains $x-\bigl((\om^i z)^b\,h(\om^i z)$ for all $i$.
Moreover, it is generated by the latter as an $\C[[x,z]]$-ideal.
Indeed, $\sum_{i=0}^{a-1}(\om^i z)^b\,h(\om^i z)=0$,
since $z^bh(z)$ does not contain $z^u$ with $u$ divisible by $a$.
Thus, $x\in (x,z^b)$, and therefore $(z^b)$ belongs to it. 
Also, all monomials in $\pm P(x,y)=\prod_{i=0}^{a-1}
(x-(\om^i z)^b\,h(\om^i z))$ belong to this ideal. 

We arrive at the isomorphism of $\C[[x,z]]$-modules:
$$
\vph: (x,z^b)/\C[[x,z]]P(x,y)\simeq \bigl(z^b\C[[z]]\bigr)^a.
$$
For any sequence of
$\r$-submodules $\vect{M}=(M_0,\ldots, M_{a-1})$ in $\C[[z]]$,
the notation is 
$$\hat{M}=\vph_x^*(\vect{M})=
\vph^*(x\vect{M})\subset (x,z^b)$$ for
the inverse image of $(xM_0,\ldots, xM_{a-1})$
in $(x,z^b)/\C[[x,z]]P(x,y)$, which is 
a $\C[[x,y]]$-module. 
The map $\vph_x^*$ is an isomorphism from
$\vect{M}$ to $\hat{M}$. The latter are called 
{\it relatively standard $\C[[x,y]]$-modules\,}
if $\C[[t]]\hat{M}=(x,z^b)/\C[[x,z]]P(x,y)$. Equivalently,
{\it all\,} $M_i$ must be standard. Note that 
$\dim (x,z^b)/\hat{M}=\sum_{i=0}^{a-1}\dim \C[[z]]/M_i$. 
\vskip 0.2cm

{\sf A simple example.}  Let $\r=\C[[x\!=\!z^3,y\!=\!z^2]]$.
Then $\de=1$, $P=x^2-y^3=(x-z^3)(x+z^3)$,\  
$(x,z^3)/\bigl(\C[[x,z]]P(x,y)\bigr)=
(x,z^3)/(x^2-z^6)$.
There are $2$ families of standard $\r$-modules:
 $M_{t\!ot}=\C[[z]]$ and $M_\al=\r (1+\al z)$ for $\al\in \C$.
They contain $\C[[z]]z^2$.
One has:\, $\vph_x^*(M_\al,M_\be)= \C[[x,y]](x+z^3)(1+\al z)+ 
\C[[x,y]](x-z^3)(1+\be z)= \C[[x,z^2]](x+(\al+\be)zx)+
\C[[x,z^2]](z^3+(\al-\be)z^4)$. This $\C[[x,z^2]]$-module contains 
$(xz^2,z^5)/(x^2-z^6)$;\, $\dim (x,z^3)/(xz^2,z^5)=4$.  
Its codimension in $(x,z^3)$ is $2$. 
 Similarly,
$\vph_x^*(M_\al,M_{t\!ot})$ and, correspondingly,
 $\vph_x^*(M_{t\!ot},M_{\al})$
are generated as $\C[[x,z^2]]$-modules
by 
$\{x+2\al zx, x-z^3,zx-z^4\}$ 
and by $\{x+2\al zx, x+z^3,zx+z^4\}$; they contain 
(the image of) $\C[[x,z]]z^5$ and are of codimension $1$ in
$(x,z^3)$. Finally, $\vph_x^*(M_{t\!ot},M_{t\!ot})=(x,z^3).$\sq
%\vskip 0.2cm

We will use the {\it full lifts} (inverse images)
 $\tilde{C}_\r$ and
$\tilde{C}^x_\r$ of $(z^{2\de})$ and $(z^{2\de+b})=x(z^{2\de})$
in $\C[[x,y]]$. Also, $\widehat{C}_\r\equal \C[[z]]\tilde{C}_\r+
\C[[x,z]]P(x,y)$ and $\widehat{C}^x_\r\equal \C[[z]]\tilde{C}^x_\r+
\C[[x,z]]P(x,y)$ are their full lifts to ideals in $\C[[x,z]]$.
Note that $\tilde{C}^x_\r \neq x\tilde{C}_\r$ and
$\hat{C}^x_\r \neq x\hat{C}_\r$.
 
\begin{proposition}\label{MIast}
We assume that the ideal 
$\tilde{C}^x_\r$ is constant for some family of rings $\r$
considered up to automorphisms 
$y\mapsto f y$, $x\mapsto g x$ for invertible 
$f\in \C[[y]]$, $g\in \C[[x,y]]$.
Then the schemes
%$J^1=\{M\mid \C[[z]]\supset M \supset \r\}$ and
$\overline{Jac}$ and $Jac^\bullet$ considered up to isomorphisms 
are constant for such a family.
%are fully determined by  $\tilde{C}_\r$, which in its turn
%determines the corresponding $\Ga$. Accordingly,  
%if there is a family of rings $\r$ 
%such that $\tilde{C}_\r$ is constant for this family
%considered up to automorphisms of $\C[[x,y]]$, then 
%the schemes $J^1, Jac^\bullet$
%are constant up to isomorphisms too. 
\end{proposition}
{\it Proof.}  
For $d\le \de$, and sequences 
$\vect{M}=(M_i\subset \C[[z]], 0\le i\le a-1)$ 
of $\r$-modules, let  $d_i\!=\!dim\, \C[[z]]/M_i$, and
$$\mathbb{J}(d)\equal \bigl\{\vect{M}\, \mid
\,\sum_{i=0}^{a-1} d_i\!=\!d \bigr\},\ \,
J_{d_0,\ldots,d_{a-1}}= \prod_{i=0}^{a-1}
J(d_i),$$
where  $J(d)\equal 
\bigr\{M\subset \C[[z]],\,
dim\,\C[[z]]/M\!=\!d\bigl\}$.

If $M_i$ are all standard,
the notation is $\mathbb{J}(d)^\bullet$, 
$J^\bullet_{d_0,\ldots,d_{a-1}}$.
 Obviously 
$\mathbb{J}(d)$ and $\mathbb{J}(d)^\bullet$ are
disjoint unions of the schemes $J_{d_0,\ldots,d_{a-1}}$ and 
$J_{d_0,\ldots, d_{a-1}}^\bullet$, where 
$\sum_{i=0}^{a-1} d_i=d$.

The map $\vph_x^*: \vect{M}\mapsto \hat{M}$ identifies the set
of all $\vect{M}$ with the set of  all
$\C[[x,y]]$-modules $\hat{M}$ in $(x,z^b)$ containing the ideal 
$\C[[x,z]]P(x,y)\!=\!\sum _{i=0}^{a-1} 
z^i\C[[x,y]]P(x,y)$\, for $P(x,y)$ from Lemma \ref{Peq}.

Given $d\le \de$,
the lift $\vph_x^*$ results in the isomorphism of schemes:
\begin{align*}
%&J(d)\equal \bigr\{\vect{M}\mid \C[[z]]\supset M_i\supset z^{2\de}\C[[z]],\,
%dim\,\C[[z]]/M\!=\!d\bigl\} 
%\simeq \\
&\mathbb{J}(d)\simeq\hat{J}(d)\equal 
\bigr\{\hat{M}\mid (x,z^b)\supset \hat{M} \supset \widehat{C}^x_\r,\ 
dim\,\C[[x,z]]/\hat{M}\!=\!d\bigl\}
\end{align*}
for all such 
$\C[[x,y]]$-modules $\hat{M}$.
We use
that $dim\,\C[[z]]/M\!\le \de$ implies
$M\supset z^{\de}\C[[z]]$; see (\ref{cjacde}). 

The assumption in the proposition
provides that $\widehat{C}^x_\r=\C[[z]]\tilde{C}^x_r$
is constant up to
isomorphisms of $\C[[x,z]]$ preserving $\C[[x,y]]$. 
Thus $\hat{J}(d)$ are topological invariants
of $\r$ and so are  $J(d)$ for any $d\le \de$.
Indeed, $\mathbb{J}(d)$ is the disjoint union
of $a$ copies of $J(d)$ and some products in terms of $J(d')$ for
$d'<d$, so we can proceed by induction
with respect to $d$ starting with $d=0$: $J(0)=\{M_{t\!ot}\}$.
We obtain that
$\overline{Jac}=J(\de)$ and, similarly, $Jac^\bullet$
are topological invariants. \sq 

There is a description of 
standard $M$ 
based on $J_1\equal\bigr\{
z^{2\de}\C[[z]] \subset{M}\subset \r\,\bigr\}$, which is the
$\ast$-dual of $J^1$ from (\ref{UJac}). Namely,
$Jac^\bullet\simeq U_1 J_1\simeq U_1 (xJ_1)$ for $U_1=1+zC[[z]]$.
Let
\begin{align}\label{Jsub1}
 \tilde{J}^x_1\equal\bigr\{\,\hbox
{ideals\,\,} \tilde{I}\subset x\C[[x,y]]\,\mid\, \tilde{I}
\supset \tilde{C}^x_\r \,\bigr\}.
\end{align}
Then $\mathbb{J}^\bullet$ is isomorphic to 
$U_1\,\tilde{J}^x_1$ generally, and  therefore
depends only on the topological type
of $\r$ under the assumption of Proposition \ref{MIast}. 
}

\comment{
Using this approach, we can obtain a somewhat different
version of Proposition \ref{MIast} for $\Jac^\bullet$,
which does not require the extension of $\C[[x,y]]$ by $z$;
its localization is used instead.

\begin{proposition}\label{MIast}
We assume that the ideal 
$\tilde{C}_\r$ is constant for some family of rings $\r$
considered up to automorphisms of $\C[[x,y]]$. Also, let
$g_i,g_{i+1}\in \tilde{C}_\r$ be some elements with 
$\nu_z=i,i+1 (i\ge 2\de)$ in $\r$ that are constant for this family.  
Then the scheme $Jac^\bullet$ considered up to isomorphisms 
is constant too.
\end{proposition}
{\it Proof.}  
One can lift the
normalization ring $\C[[z]]$ to a subring of the localization
$\C[[x,y]]_{(g_{i})}$ of $\C[[x,y]]$ 
by $g_{i}^{\Z_+}$ as follows:
$$
\tilde{O}_\r=\sum_{m=0}^{2\de-1} \,
\{\C[[x,y]](g_{2\de+1}/g_{2\de})^m\}\,+\, \tilde{C}_\r .
$$ 
It is generally not a full lift, but its image in $\r$ is the
whole $\C[[z]]$. By considering inverse images of
the standard modules $M\subset \C[[z]]$ in $\tilde{O}_\r$, 
we identify 
%We use that any elements of $\C((z))$ can be lifted to
%$\C[[x,y]]_{loc}$\,.  
$Jac^\bullet$ with the scheme
of $\C[[x,y]]$\~invariant submodules $\tilde{M}$ in
 $\tilde{O}_\r$ such that $\tilde{O}_\r= 
\tilde{M}\tilde{O}_\r$ and $P(x,y)\in \tilde{M}$. 
The latter condition potentially involves $P(x,y)$; it
can be bypassed as follows.
Using (\ref{UJac}),

\begin{align}\label{J1embed}
J^1\simeq \tilde{J}^1\equal
\bigl\{\,\C[[x,y]]\hbox{-submodules\, } \tilde{M}\subset \tilde{O}_\r
\mid \tilde{M}\supset \C[[x,y]]\,\bigr\},\\
Jac^\bullet\simeq \tilde{Jac}^\bullet\equal
\tilde{U}_\r\, \tilde{J}^1 \hbox{\, for the group $\tilde{U}_\r$ of
invertibles in\, } \tilde{O}_\r. \notag
\end{align}  
Thus $\tilde{U}_\r\, \tilde{J}^1$ and therefore
$Jac^\bullet$ considered up to isomorphisms are 
constant for the family from the proposition.

The description of $J^1$ becomes directly in
terms of $\tilde{C}_\r$ using
the duality map  $\ast$, which 
identifies $J^1$ with the scheme $J_1$ of all ideals $I\subset \r$ 
such that $ I\supset (z^{2\de})$. 
The full lifts of the elements of $J_1$ to ideals
in $\C[[x,y]]$ form the following scheme:
\begin{align}\label{Jsub1}
 \tilde{J}_1=\bigr\{\,\hbox
{ideals\,\,} \tilde{I}\subset \C[[x,y]]\,\mid\, \tilde{I}
\supset \tilde{C}_\r \,\bigr\}\simeq J^1.
\end{align}
Here \,$\tilde{J}_1\subset \cup_{n\le \de\,} H^{(n)}$, so  
it is a disjoint union of schemes 
corresponding to different $n$ here. Then $Jac^\bullet$
becomes isomorphic to $\tilde{U}_\r \,\tilde{J}_1$.
}

\vskip 0.2cm
{\sf Semigroups with $2,3$ generators.} The assumption of the 
proposition is sufficiently explicit for such semigroups
$\Ga$. First of all,
let $\Ga^{(2)}=\lan a,b\ran$ for $1<a<b\in \Z$ such that $gcd(a,b)=1$.
It is for $\r^{(2)}=\C[[y=z^a, x=z^b h(z)]]$ and any 
$h(z)=1+z(\ldots)$
as above. By $\lan \cdots\ran$, we mean here the span over $\Z_+$.

Now we consider the case of
$3$ generators. Let $1<a\!=\!\upsilon m<b\!=\!\upsilon n$ for 
$\upsilon\!=\!gcd(a,b)>1$,
$\r=\C[[y\!=\!z^{\upsilon m}, x\!=\!z^{\upsilon n}h(z)]]$, where
$h(z)=1+Cz^{p}+\mod (z^{p+1})$ for $\C\neq 0$, and $b+p$ 
is the first $z$-exponent
in $z^b h(z)$ that is not a linear combination over $\Z_+$ of $a,b$.
Then $x^m-y^n=mCz^{\upsilon mn+p} 
\mod (z^{\upsilon mn+p+1})$ and 
therefore $\upsilon mn+p\in \Ga$ is the smallest element not
in the $Z_+$-span of $a,b$. 

Let $gcd(p,\upsilon)=1$, which is necessary and sufficient
to ensure that 
 $\Ga$ has $3$ generators: they 
will be $\upsilon m, \upsilon n, 
\upsilon mn+p$.
We will denote such ring $\r$ and the semigroup 
$\Ga=\lan \upsilon m, \upsilon n, \upsilon mn+p\ran$ by 
$\r^{(3)},\Ga^{(3)}$. In this case, $\upsilon n+p$ 
is the {\it Zariski invariant} of $\r$.  Here
we can (and will) assume that $h(z)=1+z^{p}+uz^{p+\kap}\mod 
(z^{p+\kap+1})$
for some $u$, where $\kap\!\equal\!
Min[(\Z_+\!\setminus\! \La')\cap (p\!+\!\Z_+)]\!-\!p$
for $\La'$ from (\ref{Laprime}). We set $\kap=0$
if this intersection is empty; then $\Ga$ fully determines the
singularity.

The formulas for the corresponding $2\de$ are as follows:
\begin{align}\label{de12}
2\de^{(2)}=(a\!-\!1)(b\!-\!1),\, \,2\de^{(3)}=
\upsilon^2 mn\! -\!\upsilon(n\!+\!m)\!+\!(\upsilon\!-\!1)p\!+\!1.
\end{align}
We use here formula (2.1) from \cite{HH}
for $2\de$, which is for any $\Ga$. By the way, the latter gives
that the condition $gcd(p,\upsilon)=1$ implies that $\Ga$
has exactly $3$ generators. Otherwise the corresponding $2\de$
would be smaller than $2\de^{(3)}$ above, and other 
generators must be no smaller than
$\upsilon(\upsilon mn+p)>2\de^{(3)}$, which is impossible
(they cannot belong to $\Ga^{(3)}$).

See \cite{HH} for the (formal)
analytic classification of the singularities
with $a=4,b=9$ and the table for $a=4$. For instance, 
$\r_u=\C[[y\!=\!z^4, x\!=\!z^9\!+\!z^{10}\!+\!uz^{11}$
are non-isomorphic for different $u\in \C$. Here $10$
is the {\it Zariski invariant} and 
$\La'=\nu_z\bigl(\r z^{a-1}+\r (dx/dz)
\bigr)+1-a=\lan 0,4,5,8,9,10,12,13,14,15,\ldots\ran$
for $u\neq \frac{19}{18}$.
We use that $f\!\equal\!z^4(dx/dz)\!-\!9xz^3\!=\!z^{17}\!+\!2u z^{18}$ 
and $x(dx/dz)\!-\!9yf=(19\!-\!18u)z^{18} \mod (z^{19})$.
Only $11$ is missing after $10$, so adding the terms 
$z^{> 11}$ to $x$
does not change the corresponding analytic type of $\r$.  
For $\upsilon>1$, an example of $\r$ with
continuous moduli is  $a=4,b=10,p=1$. 
\vskip 0.1cm

Let 
Ceiling$\,[e]=\min\{\Z\ni n\ge e\}$ in the next theorem.
The polynomial $P$ is as above; it depends only on $a,b,p$ 
for any $\r^{(3)}$ with $\kap=0$.

\begin{theorem}\label{thmtop}
(i) 
The polynomial $P(x,y)$ is constant modulo the monomials 
with the $\nu_z$-valuations no smaller than $2\de$  
for any $\r^{(2)}$ and for
$\r^{(3)}$ with $\kap\ge 1$ under the following inequalities 
(for $p$):   
\begin{align}\label{ineqp}
&(\upsilon-1)p\le a\cdot\!
\hbox{Ceiling\,}[\,\frac{b+p+\kap}{a}\,]+a\!-\!1, \text{\ or}\\
%(\upsilon-1)p<\upsilon(n+m)=a+b \hbox {\, in the case of\,}\, \r^{(3)},
&\hbox{for a stronger one\,}\,\,(\upsilon-2)p\,\le\, a+b+\kap\!-\!1\,. 
\label{ineqp1}
\end{align}
This gives that
 $\tilde{C}_\r$ coincides with that for $y=z^a,x=z^b+z^{b+p}$,
Proposition \ref{MIast} is applicable,
and $\tilde{Jac^\bullet}$ is a topological invariant of $\r$.

(ii) If $(b+p+\kap)$ is replaced by $(b+p)$ in (\ref{ineqp})  and
$a+b+\kap-1$ by $a+b-1$ in (\ref{ineqp1}), then  
$\tilde{C}_\r$ has a basis of 
eigenvectors for the action $x\mapsto v^bx,\, y\mapsto v^ay,$ 
$v\in \C^\ast$. Here the case $\kap=0$ (when $\r$ is
determined by $\Ga_\r$) is included.
Furthermore, all monomials
of $P$ belong to $\tilde{C}_\r$ if\, $(\upsilon\!-\!1)p\le 
a\!+\!b\!-\!1$. Finally, if $p=1$ for $\r^{(3)}$ or for
any $\r^{(2)}$, then the ideal $\tilde{C}_\r$ is monomial (and depends
only on $\Ga_\r$).
\comment{
Assuming that 
$(\upsilon-1)p-a\,\hbox{Ceiling\,}[\frac{b+p+1}{a}]\le a-2$, 
e.g. $(\upsilon-2)p\le a-1$, 
$\tilde{C}^x_\r=x\tilde{C}_\r+ (x^a-y^b)\C[[x,y]]$, and
$\tilde{C}_\r,\tilde{C}^x_\r$ depend only on $\Ga^{(3)}$, i.e. 
do not depend on $h(z)$. Thus 
Proposition \ref{MIast} is applicable,
and $\overline{Jac}$,  $Jac^\bullet$ are topological invariants of
the corresponding singularities.
%generated by $x^a,y^b$. 
Furthermore,  $\tilde{C}_\r,\, \tilde{C}^x_\r$
are invariant under the action $x\mapsto u^bx,\, y\mapsto u^ay$
for $u\in \C^\ast$ if $(\upsilon-2)p\le a-1$, and all monomials
of $P$ belongs to $\tilde{C}^x_\r$ if $(\upsilon-1)p\le a-1$
When $p=1$ for $\r^{(3)}$ and for $\r^{(2)}$,
the corresponding   $\tilde{C}_\r$ and $\tilde{C}^x_\r$
are monomial ideals.
}
\end{theorem}
{\it Proof.} {\sf (i)} Let us examine 
the monomials in $P(x,y)$ with potentially non-constant
coefficients, i.e. those depending on $u_j (j\ge p+\kap)$ in
the presentation $x=z^b+z^{b+p}+u_{p+\kap}z^{b+p+\kap}+\,\ldots$\ . 
A natural lower bound of the $\nu_z$-valuations for such 
monomials in $P(x,y)$ is $a(b+\ep)$, where $\ep\ge \hbox{Ceiling}
\,[\,\frac{b+p+\kap}{a}\,]-\frac{b}{a}\ge \frac{p+\kap}{a}$; we use
Lemma \ref{Peq}. This readily gives 
that (\ref{ineqp})
or its somewhat stronger version (\ref{ineqp1}) provide that
$P(x,y)$ is constant modulo $\tilde{C}_\r$.
%The deduction of the first inequality from the second one is
%due to the following obvious bound:
%$\hbox{Ceiling\,}[\,\frac{b+p+1}{a}]\ge \frac{b+p+1}{a}$.
Thus, $P(x,y)$ coincides with  $P_0(x,y)$ 
obtained for $x=z^b+z^{b+p}$ 
modulo the monomials with $\nu_z\ge 2\de$. 

To conclude $(i)$, $\tilde{C}_\r$ is constant,
since we can find at least one {\it constant\,}
element $g_k\in \C[[x,y]]$ for any $k\ge 2\de$ with 
$\nu_z(g_k)=k$.  These elements will be in $\tilde{C}_\r$
by the definition of the latter.
Namely, one represents:
$k=\al a+\be b+\ga(\upsilon mn+p)$ with $\al,\be,\ga\in \Z_+$
for any $k\ge 2\de$ ($\ga=0$ for $\r^{(2)}$). Then we set
$g_k=y^\al x^\be (x^m-y^n)^\ga$.
Together with the ideal $\C[[x,y]]P_0(x,y)$, such 
elements $g_k$ for $k\ge 2\de$
generate $\tilde{C}_\r$. 
\vskip 0.2cm

{\sf (ii)}  First of all,
the inequality $2\de\le ab$, which provides that all
monomials of $P$ belong to $\tilde{C}_\r$, is equivalent
to $(\upsilon-1)p\le a+b-1$ from $(ii)$.
We use that $2\de^{(3)}=(a-1)(b-1) +(\upsilon-1)p$. 
Technically, we can set $p=1,\upsilon=0$ 
for $\r^{(2)}$; obviously, $2\de=(a-1)(b-1)< ab$ in this case.

Next, the first
inequality from part $(ii)$ gives that $\tilde{C}_\r$
is generated by the elements of valuations no smaller than $2\de$
and the monomials
of  $P$ of the valuation $ab$ (only such are sufficient). 
This gives that $g_k$ above 
are eigenvectors under the action from $(ii)$.

\vskip 0.2cm
 
{\sf The monomiality for $p=1$ and $\r^{(2)}$.}
Since $\ga$ used above (for $g_k$) is $0$ for $\r^{(2)}$, 
$\tilde{C}_\r$ is a monomial
ideal in this case. Namely,
it contains $x^a,y^b$, which have the valuation $ab$, and
$g_k$ are monomials for $(a\!-\!1)(b\!-\!1)\le k<ab $, 
which is sufficient. We note that they are unique 
monomials with  $\nu_z=k$ up to proportionality
in this range of $k$ .

%\vfil
Let us consider $\r^{(3)}$ with $p=1$.
One has: $2\de/\upsilon=\upsilon mn -m-n+1$ for $p=1$.
Let $(m-1)(n-1)=\al m+\be n$ with $\al,\be\in \Z_+$.
We obtain that
$2\de/\upsilon=(\al+(\upsilon-1)n)m+\be n=
 \al m+ (\be+(\upsilon-1)m)n.$ Therefore,
$g_y=y^{\al+(\upsilon-1)n}x^\be$ and  
$g_x=y^\al x^{\be+(\upsilon-1)m}$; both elements are
with $\nu_z(g)=2\de$.
One has:
$$
g_x\!-\!g_y=y^\al x^\be (x^{(\upsilon-1)m}\!-\!y^{(\upsilon-1)n})=
x^\be y^{\al+(\upsilon-1)n}\bigl((\upsilon\!-\!1)m z\!+\! z^2(\cdots)
\bigr).
$$
Thus, $g_x-g_y$ represents $\nu_z=2\de+1$, where the monomials
$g_x$ and $g_y$  are
in $\tilde{C}_\r$. (This step is more involved 
in the nonzero characteristic.) 

We note $(\upsilon-1)mn$ is the smallest number such that it has
$\upsilon$ different $\Z_+$\~representations in terms of $m,n$
according to Theorem 4 from
\cite{BR}. It gives that some 
linear combinations of the monomials in $\tilde{C}_\r$
provide the  
valuations $2\de+1,\ldots, 2\de+(\upsilon-1)$, which is sufficient
to finish the proof. 
However, a simpler argument works.

We represent  $2\de\!+\!q$ with $q\ge 0$ as 
$\al' (\upsilon m)\!+\!
\be'(\upsilon n)\!+\! \ga' (\upsilon mn+1)$ with 
$\ga'=q \mod \upsilon$. The last condition is necessary here.
It is sufficient because $\upsilon (\upsilon mn+1)$ can
be represented in terms of $a,b$ over $\Z_+$; any number
no smaller than $\upsilon (m-1)(n-1)$ is such.  Thus, $\nu_z$
of any monomial in the product  $y^{\al'}x^{\be'}
(x^m-y^n)^{\ga'}$ is $\al' (\upsilon m)+
\be'(\upsilon n)+ \ga' (\upsilon mn)=(2\de+q)-\ga'\ge 2\de$. So
they  belong to the monomial part of $\tilde{C}_\r$. \sq

\vskip 0.2cm
{\sf An example of $p=1$.}
Let us illustrate the theorem for  $a\!=\!6,b\!=\!9,p\!=\!1$, i.e.
for $\upsilon\!=\!3, m\!=\!2, n\!=\!3$ and 
$\r=\C[[y\!=\!z^6, x\!=\!z^9h(z)]]$. One can assume
that $h(z)=1+z+\ldots$.  
Then $\nu_z(x^2-y^3)=
19,\,  2\de=42$, where
\begin{align*}
\tilde{C}_\r=\bigl\langle &\{y^7,x^4y,x^2y^4\}_{42}, 
\{y(x^4-y^6), y^4(x^2-y^3)\}_{43},
\{y(x^2-y^3)^2\}_{44},\\
&\{x^5, x^3y^3, xy^6\}_{45}, \{x^3(x^2-y^3)\}_{46},
\{x(x^2-y^3)^2\}_{47}\ldots\bigr\rangle;
\end{align*}
we show only some elements for the corresponding
 $\nu=\nu_z$ (inside  $\{\,\}_{\nu}$). Therefore:
$\tilde{C}_\r=\lan y^7,x^4y,x^2y^4, x^5, x^3y^3, xy^6\ran$ as
a $\C[[x,y]]$-module.
The corresponding $\la$ and $\la'$ for this monomial ideal are
\Yboxdim5pt
$$
\la(\tilde{C}_\r)=\{7,6,4,3,1\}, \la'= \yng(7,6,4,3,1). 
$$
\Yboxdim7pt
The number of boxes is $\de=21$. 
For $x=z^9+z^{10}$, one has:\, 
$P=y^{10}-6xy^8-y^9-2x^3y^5+3x^2y^6-3x^4y^3+x^6$, where
all monomials have their valuations no smaller than
$54$; so they are well inside $\tilde{C}_\r$.

{\sf An example of $p=2$.}
The ring $\r=\C[[y\!=\!z^6, x\!=\!z^9+z^{11}]]$ 
is the simplest when some
Piontkowski cells are non-affine spaces; see the Appendix to
\cite{ChP1}. In this case, $\Ga=\lan 6,9,20=\upsilon n+p\ran$, 
$2\de=(a-1)(b-1)+(\upsilon-1)p=44$,
$ab=54>2\de=44$, and 
$a+b-1=14\ge p(\upsilon-1)=4$. The latter inequality gives
that all monomials of $P$ belong
to $\tilde{C}_\r$.
Let us provide $P(x,y)$ for this $\r$:
$
x^6 - 3 x^4 y^3 + 3 x^2 y^6 - 6 x^2 y^7 - y^9 - 2 y^{10} - y^{11}.
$
In this case, $\tilde{C}_\r$ is generated by
$\{(x^2-y^3)x^2y\}_{44}, \{x^5\}_{45}$,
$\{(x^2-y^3)^2y\}_{46}$,
$\{(x^2-y^3)x^3\}_{47}, \{y^8\}_{48}, \{(x^2-y^3)^2x\}_{49},
\ldots\ $; we
show the corresponding valuations.
To obtain the other representatives, multiply those we provided
by powers of $y$. This is not a monomial
ideal.
\vskip 0.2cm

Recall that always $P(x,y)\in \tilde{C}_\r$ by definition,
and that this ideal generally
depends on the coefficients of $P$ if $2\de> ab$.
Even if all monomials of $P$
belong to $\tilde{C}_\r$, its linear generators
can involve the coefficients of $h(z)$. On the other hand,
if Theorem \ref{thmtop} is not applicable, the 
 ideal $\tilde{C}_\r$ can be still a topological invariant
of $\r$ up to isomorphisms, as well as  
$\tilde{Jac^\bullet}$.
A significantly weaker propriety
was conjectured in \cite{ChP1}:\  that   
the motivic super-polynomials are always topological invariants
(for plane curve singularities), 
i.e. depend only on $\Ga_\r$.
\vskip 0.1cm 

We think that this theorem provides the main cases when 
one can obtain a reasonably simple connection between
the refined invariants of plane curve singularities 
from \cite{ChD1} and Hilbert schemes of $\C^2$.
We mean here mostly some possible
generalizations of Theorem 1.1 from \cite{GoN}. 
It is not surprising
that their theorem was restricted to torus knots; in this case,  
$\tilde{C}_\r$ is the simplest.

\subsection{Employing the parametrization}
Given any $M\in \overline{Jac}$, 
we set $M^e\equal z^{2\de+e}M$, which is an ideal in $\r$
containing $(z^{4\de+e})$ for $e\in \Z_+$.
Recall that the conductor $c(M)$ is no greater than 
$\de+|Z_+\setminus \De(M)|$ for standard $M=M_\bullet$, i.e. 
$M\supset (z^{\de+|Z_+\setminus \De(M)|})$;
see (\ref{rval}) and \cite{PS}. Since 
$\overline{Jac}=\{z^{|\De(M)\setminus \Ga|}M\}$ for
standard $M$, the modules for the points of  $\overline{Jac}$
automatically contain $(z^{2\de})$. 

\comment{ 
This gives that $M\in \overline{Jac}$, which can be represented
as $z^{|\De(M)\setminus \Ga|}M$ 
$$
M=z^{}M_\bullet, \where dev(M_\bullet)\equal
|\De_\bullet\setminus \Ga|=\de-|Z_+\setminus \De_\bullet|,
$$
and $M_\bullet\supset (z^{\de+|Z_+\setminus \De_\bullet|})$;
see \cite{PS}.
}

The embedding   $M^e\subset \r$ 
can happen for $e<0$ for some $M$, which values 
of $e$ will be allowed in the considerations below. For such $e$:
$$\dim_{\C} \r /M^e\ =\ 
\dim_{\C} \C[[z]]/M^e\!-\!\de\ =\ 2\de+e.$$

The (full) lift of $M^e$ to the ideals from $\C[[x,y]]$ is 
natural:
\begin{align}\label{MtoI}
\tilde{I}^e(M) \equal \C[[x,y]]P(x,y)+\C[[x,y]]\tilde{M'},\ 
\tilde{I}(M)\equal \tilde{I}^{e=0}(M)
\end{align} 
for any set $\tilde{M}'\subset \C[[x,y]]$ such that its image in $\r$
is $M^e.$ One has:
$\dim_{\C}\C[[x,y]]/\tilde{I}^e(M)=2\de+e$, and 
$\tilde{C}_\r=\tilde{I}(\C[[z]])$ in this notation.
We will use that
$x^i y^j\in \tilde{I}^e{M}$\, if\, $bi+aj\ge 4\de+e$. 

Using the operation $I\mapsto I_0$ for the  
ideals $I\subset\C[[x,y]]$
from the previous section, we set: 
$$
\tilde{I}^e_0(M)
\equal (\tilde{I}^e(M))_0=I_{\tilde{\la}^e}
$$
for the corresponding partition $\tilde{\la}^e=
\tilde{\la}^e(M)$ of
order  $2\de+e$. Since $P(x,y)\in \tilde{I}^e(M)$, it
has no greater than $a$ lines, i.e. 
$\ell(\tilde{\la})\le a$. Given any $\r$\~module $M$,
the polynomial $P(x,y)$ is uniquely determined by 
$\tilde{I}^e(M)$ for sufficiently large $e$, which is 
part of the following proposition.

\begin{proposition}\label{PfromI}
(i) Given $M\in \overline{Jac}$, the polynomial $P(x,y)$
coincides with $f_{a0}$ of $\tilde{I}^e(M)$ for sufficiently
large $e$. Combinatorially, this holds if and only if
the diagram $\tilde{\la}'$ for $\tilde{\la}^{e}(M)$ 
contains all monomials from $P(x,y)-(-x)^a$, i.e. all boxes
$\{i,j\}$
with $d_{i,j}\neq 0$ in  (\ref{Peq}). 

(ii) For such $e$, $\tilde{I}^{e+a}(M)=y \tilde{I}^{e}(M) 
+C[[x,y]]f_{a0}$, and the diagram for $\tilde{\la}^{e+a}(M)$ is 
obtained from $\tilde{\la}'$ for $e$  by adding one 
$a$\~column at its beginning. In the presentation from Theorem 
\ref{thmc}, the polynomials $p_{i,j}$ remain unchanged when
$e\mapsto e+a$,  but the corresponding 
$y^{d_1}$ must be replaced for $e+a$ by $y^{d_1+1}$.   \sq 
\end{proposition}

The ideals $\,\tilde{C}^{(e)}_\r\equal \tilde{I}^e(\C[[z]])\,$ 
are the key in 
the description of the $\overline{Jac}\,$ via $\,C[[x,y]]$, 
including $\tilde{C}_\r^{(0)}=\tilde{C}_\r\,$ used in
Proposition  \ref{MIast}. 
Note that $e=0$ is generally far from the "stabile
values" of Proposition \ref{PfromI}. However it
is of importance to calculate the $\la$-partitions  for 
$e=0$ and for 
$\tilde{I}^{min}_M$, which is defined as $\tilde{I}^{e}(M)$ with the 
smallest possible  $e$ for one or some $M$
providing that $z^{2\de+e}M\subset \r$. The
Piontkowski strata $\{M\!=\!M_\bullet \mid \De(M)\!=\!\De\}$
are natural here to analyze.

For an individual $M$, 
it is actually more natural to allow here more relaxed
embeddings
$\phi_{2\de+e'} M\subset
\r$ for {\em proper\,} 
$\phi_{2\de+e'}\in z^{2\de +e'}+ (z^{2\de+e'+1})$.
According to (\ref{u-min}), one has:
%$2\de+e_\circ'=\min \De(M^\ast)= 2\de -c(M)$ and
$\min \{e'\}=-c(M)$ 
for the conductor $c(M)$. This relaxation makes sense
for families of $M$
too, assuming that  $Aut(M)$ is fixed in this family; cf. 
\cite{Sto}. 
\vskip 0.2cm

Let us summarize what we obtained, 
combining our analysis with Theorem \ref{thmc}.

\begin{theorem}\label{thmtoI}
For a ring $\r\subset \C[[z]]$ with the generators $x,y$
picked as in Lemma \ref{Peq}, let  $\Ga$, $\de$, and the
compactified Jacobian $\overline{Jac}\ni M$ be as above.
Let $\pi^e(M)=\tilde{I}^e(M)$, which is the inverse image
of $M^e=z^{2\de+e}M$ in $\C[[x,y]]$, an ideal 
of codimension $2\de +e$.  
Here $e\ge 0$ ensures that $M^e\subset \r$ for all
$M\in \overline{Jac}$, but it can be
negative for special families of $M$. 
As above:\ 
$\tilde{C}^{(e)}_\r=\tilde{I}^e(\C[[z]]),\ \tilde{C}_\r=
\tilde{C}^{(0)}_\r$.

(i) The map $\pi^0$ establishes an isomorphism
of schemes between $\overline{Jac}$ and
$\pi^0(\overline{Jac})\subset H^{(2\de)}$.
The latter coincides with the projective subscheme
$\{I\in H^{(2\de)}\mid P(x,y)\in I\subset \tilde{C}_\r\}$ of
$H^{(2\de)}$ for the equation $P(x,y)$ as above.
Such ideals $I$ automatically contain
$\tilde{C}^{(2\de)}_\r$, which is  
$(z^{4\de})$ lifted to $\C[[x,y]]$, including 
all $x^iy^j$ such that $bi+aj\ge 4\de$. Recall that 
 $H^{(2\de)}$ contains
all monomials from $\mathfrak{m}^{2\de}$, i.e. for $i+j\ge 2\de$.

(ii) Let us fix $\la$ such that $|\la|=2\de+e$ and $\ell(\la)=a$,
where $\,e$ is assumed to satisfy the stabilization conditions from
Proposition \ref{PfromI}.
Then $I\in C_\la$ for the Gr\"obner cell $\,C_\la\subset
H^{(2\de+e)}$ can be presented as $\pi^e(M)=\tilde{I}^e(M)$ for 
some $M\in \overline{Jac}\,$ if and only if 
$$ (a)\ I\subset \tilde{C}^{(e)}_\r=\tilde{I}^e(\C[[z]]) 
\and (b) \ f_{a0}=T_\la^{(1)}=P(x,y).$$ 
Accordingly, 
$C_\la\cap \pi^e(\overline{Jac})$ is the scheme of common
zeros of all polynomials $\{p_{i,j}\}$
from Theorem \ref{thmc}.

(iii) Continuing,
the generators $f_{i^\circ\!j^\circ}$ for
the corners $\{i^\circ,j^\circ\}$ of $\la'$ with $j^\circ>0$ 
are as follows. We switch to $\hat{\la}$  obtained
from $\la$ by removing the first column in $\la'$
and consider the 
corresponding generators $\hat{f}$, given by the minors 
$T_{\hat{\la}}$
with the corresponding  $p$\~polynomials obtained from 
$\{p\}$ as follows: $\hat{p}_{i,j}=p_{i,j+1}$, where $j\ge 0$. 
Then $\{i^\circ, j^\circ\! -\!1\}$ constitute  the set of
 corners of $\hat{\la}'$ and 
$f_{i^\circ\!j^\circ}=y 
\hat{f}_{i^\circ\!(j^\circ\!-\!1)}$
for $j^\circ>0$. 
\sq
\end{theorem}

\vskip 0.2cm

\subsection{Quasi-homogeneous singularities}
They are for the rings \,\, $\r=\C[[x\!=\!z^r,y\!=\!z^s]]$ for $r>s>0$ 
such that
$gcd(r,s)=1$. Then $\de=(r-1)(s-1)/2$ and the lift
of $(z^{2\de})$ to $\,\C[[x,y]]$ is 
$$\tilde{C}=\tilde{C}_\r=
\{\sum_{i,j}c_{ij} x^i y^j \mid i,j\in \Z_+, ir+js\ge 2\de\}.
$$
This is a monomial ideal, which dramatically simplifies
the usage of punctual Hilbert schemes for $\C[[x,y]]$ for
the study of $\overline{Jac}$ and $Jac^\bullet$.

The diagram $\la'$ for the
partition $\la$ of $\tilde{C}$
is  formed by {\em all\,}
boxes in the $s\times r$\~rectangle 
below the diagonal connecting $\{i=s\!-\!1,j=0\}$ and
$\{i=0,j=r\!-\!1\}$ in our standard presentation
of diagrams.  Their number is
indeed $(r-1)(s-1)/2$. As above, the boxes are numbered by their
upper-left corners $\{i,j\}$, where  $\ 0\le i\le s\!-\!1,\ \,
0\le j \le r\!-\!1$.

Let us put the numbers $(2\de-1)-ir-js$ in the corresponding boxes.
See Figure \ref{3x4-diag} for $r=4,s=3$. Then we arrive at
the interpretation of the Piontkowski $\De$-modules from 
\cite{Pi} in terms of the {\em Dyke paths\,} from \cite{GM1}; 
see there Section 2.2 and Figure 1. Let us state it and
connect it with  $\tilde{\la}^{min}(M)$ 
for {\em $z$\~monomial\,} 
$\r$\~modules $M\subset \C[[z]]$.

Recall that $\De_\bullet(M)=
\De(M_\bullet)=\De(M)-v$ \,for $M_\bullet=z^{-v}M$, 
where $v=\min{\De(M)}$, where $\De(M)$ is the
$\Ga$\~modules of $M$: \,
$\De\subset \Z$ and  $\Ga+\De\subset \De$. 
We always assume that $\De\subset \Z_+$. 
Any $\Ga$\~modules in the
quasi-homogeneous case come from some $M\in \C[[z]]$.
They are fully described
by their sets of  {\em gaps\,}, which is 
$\Z_+\setminus \De$. Finally, $\hat{\De}\equal
\De_\bullet\setminus \Ga$, which is the set of {\em added gaps}
 from $\Ga$ for 
$\De_\bullet=\De-\min{\De}$.  
\vskip 0.2cm

We define {\em Dyck paths\,} 
as Young diagrams in this rectangle, which can be empty,
above the anti-diagonal, i.e. with $ir+sj<2\de$ 
in the $\{i,j\}$-presentation.  
The correspondence
from \cite{GM1} is between the Dyck paths 
and  standard $\De$. Namely,
the set of numbers $2\de-1-ir-js$\, calculated for 
the boxes of the corresponding Dyck path
is $\hat{\De}$. The next
proposition follows from this interpretation.

\begin{proposition}\label{Pio}
Given a $z$\~monomial standard $\r$\~module $M$,
let  $\la'$ be the diagram constructed from the partition 
$\la=\tilde{\la}^{min}(M)$. 
Then it coincides with the Dyck path for the 
(standardization of the) dual module
$M^\ast$, which is 
$\{f\in \r \mid fM\subset \r)$. \sq
\end{proposition}
\vfil 

{\sf Example: $r=4,s=3$.} Then $\r=\C[[x\!=\!z^4,y\!=\!z^3]]$, 
$2\de=6$. We will calculate $\la'$ from Proposition \ref{Pio} 
for all {\it $z$-monomial\,} standard $M$ and the corresponding 
diagrams 
$\tilde{\la}'= \tilde{\la}'(z^{2\de}M)$. Generally, the
corresponding inverse images are not monomial
in $\C[[x,y]]$ for $z$-monomial modules $M$; this is due
to the presence of $\C[[x,y]](x^a-y^b)$ in the lifts.
However they are
monomial in this particular case. 
\vfil

 First, 
$z^6\C[[z]]\subset \r$ is the minimal embedding of  
$M_{t\!ot}\equal \C[[z]]$ into $\r$. So 
$\hat{\De}_{t\!ot}=\{1,2,5\}=\Z_+\setminus \Ga$, and 
$\la_{t\!ot}'=\,\yng(2,1)\,=\tilde{\la}_{t\!ot}'$. Second, 
$\la_{0}'=\emptyset$, $\tilde{\la}_{0}'= \yng(2,2,2)\,$ for 
$M_0\equal \r$. The modules $M_{t\!ot}$,$M_0$ 
are standard selfdual. %We have 3 more cases. 
\vfil

Let $M_{1}$ be generated by $1,z$ over $\r$. It is linearly 
generated by $\{1,z,z^3,z^4,z^5,z^6\ldots\}$, so 
$\hat{\De}(M_1)=\{1,5\}.$ Its standard dual is $M_3$ below. Here 
minimal $v$ such that $z^v M_1\subset \r$ is $3$. So $z^3 M_1$ 
is the linear span $\{z^3, z^4, z^6,z^{7},\ldots\}= 
\{y,x,y^2,yx,\ldots\}$, and  %\Yboxdim7pt
$\la_1'= \yng(1)\,$, $\tilde{\la}_1'= \yng(2,1,1)$\,. The latter
is for $y\{y,x,y^2,yx,\ldots\}$, which results in the additional
$3$-column.
\vfil
 
The module $M_2\equal\lan 1,z^2\ran$ is the linear span of 
$\{1,z^2,z^3,z^4,z^5,z^6\ldots\}$.  This module 
is standard self-dual. One has: 
$\hat{\De}(M_2)=\{2,5\}$, $v=4$ and  $z^4M_2$ is the 
linear span of $\{z^4, z^6, z^7, z^{8},\ldots\}=
\{z^4, z^6, z^7, z^{8},\ldots\}$; for $z^6M_2$ it is 
$\{y^2,x^2,y^3,y^2x,\ldots\}$. Thus $\la_2'=\yng(2)$\,
and $\tilde{\la}_2'=\yng(3,2)\,$.
\vfil

Finally, let $M_3\equal \lan 1,z^5\ran$; this
conclude the list. It is the span of
$\{1,z^3,z^4,z^5,z^6\ldots\}$. It coincides with the standard
dual of $M_1$, which is
 $(M_1^\ast)_\bullet$.  Then
$\hat{\De}(M_3)=\{5\}$, 
$v=3$ and  $z^3 M_3$ is the linear span 
of $\{z^3, z^6, z^7, z^{8},\ldots\}=
\{y,y^2,yx,x^2,\ldots\}$. So $\la_3'=\yng(1,1)$\,,
 $\tilde{\la}_3'=\yng(2,2,1)$\,.
\vfil
 
It is of interest to interpret combinatorially 
$\tilde{I}^e(M)$
for {\it any\,} admissible $e$. We do this 
in this section only for minimal $e$ and $e=0$. 
Generally non-trivial
stratifications of $\overline{Jac}$ can be obtained
this way,  which are in a sense
"orthogonal" to the  
Piontkowski one in terms of $\Ga$-modules $\De$. The following
stratification of $Jac$, which is the Piontkowski stratum 
of $\overline{Jac}$ corresponding to $\De=\Ga$, is of 
particular interest. 

We consider the  
invertible $\r$-modules belonging to various 
standard {\it $z$-monomial} modules. The stratification
of $\overline{Jac}$ by the set-theoretical
differences of the corresponding {\it closures} of
these sets in $\overline{Jac}$ is presumably related to
the $q\leftrightarrow  t^{-1}$ duality of the motivic
superpolynomials.

\begin{figure}
{\noindent
\thicklines
{\makebox(20,20){}}
{\makebox(20,20){}}
{\makebox(20,20){}}
{\makebox(20,20){}}
{\makebox(20,20){j=0}}
{\makebox(20,20){j=1}}
{\makebox(20,20){j=2}}
{\makebox(20,20){j=3}}
%{\makebox(20,20){j=4}}
{\makebox(20,20){}}
{\makebox(20,20){}}
{\makebox(20,20){}}
{\makebox(20,20){}}\\
\thicklines
{\makebox(20,20){}}
{\makebox(20,20){}}
{\makebox(20,20){i=0}}
{\framebox(20,20){$1$}}
{\framebox(20,20){$y$}}
{\makebox(20,20){$y^2$}}
{\makebox(20,20){$y^3$}}
%{\makebox(20,20){$y^4$}}
{\makebox(20,20){}}
{\makebox(20,20){}}
{\makebox(20,20){}}\\
\thicklines
{\makebox(20,20){}}
{\makebox(20,20){}}
{\makebox(20,20){i=1}}
{\framebox(20,20){$x$}}
{\makebox(20,20){$xy$}}
%{\makebox(20,20){}}
{\makebox(20,20){$xy^2$}}
{\makebox(20,20){$\cdot$}}
{\makebox(20,20){}}
{\makebox(20,20){}}
{\makebox(20,20){}}\\
\thicklines
{\makebox(20,20){}}
{\makebox(20,20){i=2}}
{\makebox(20,20){$x^2$}}
%{\makebox(20,20){}}
{\makebox(20,20){$x^2y$}}
{\makebox(20,20){$\cdot$}}
{\makebox(20,20){$\cdot$}}
{\makebox(20,20){}}
{\makebox(20,20){}}

{\noindent
\thicklines
{\makebox(20,20){}}
{\makebox(20,20){}}
{\makebox(20,20){}}
{\makebox(20,20){}}
{\makebox(20,20){j=0}}
{\makebox(20,20){j=1}}
{\makebox(20,20){j=2}}
{\makebox(20,20){j=3}}
%{\makebox(20,20){j=4}}
{\makebox(20,20){}}
{\makebox(20,20){}}
{\makebox(20,20){}}
{\makebox(20,20){}}\\
\thicklines
{\makebox(20,20){}}
{\makebox(20,20){}}
{\makebox(20,20){i=0}}
{\framebox(20,20){$5$}}
{\framebox(20,20){$2$}}
{\makebox(20,20){$-1$}}
{\makebox(20,20){$-4$}}
%{\makebox(20,20){$y^4$}}
{\makebox(20,20){}}
{\makebox(20,20){}}
{\makebox(20,20){}}\\
\thicklines
{\makebox(20,20){}}
{\makebox(20,20){}}
{\makebox(20,20){i=1}}
{\framebox(20,20){$1$}}
{\makebox(20,20){$-2$}}
%{\makebox(20,20){}}
{\makebox(20,20){$-5$}}
{\makebox(20,20){$\cdot$}}
{\makebox(20,20){}}
{\makebox(20,20){}}
{\makebox(20,20){}}\\
\thicklines
{\makebox(20,20){}}
{\makebox(20,20){i=2}}
{\makebox(20,20){$-3$}}
%{\makebox(20,20){}}
{\makebox(20,20){$-6$}}
{\makebox(20,20){$\ \cdot$}}
{\makebox(20,20){$\ \cdot$}}
{\makebox(20,20){}}
{\makebox(20,20){}}
}
\vskip -0.4cm
\caption{$ \{i,j\}\mapsto 2\de-1-ri-sj$}
\label{3x4-diag}
}
\vskip -0.3cm
\end{figure}

\vskip 0.2cm
{\sf The diagrams for
$\,\overline{Jac}\,$}. Continuing with the same  $r=4, s=3$,
let us provide a complete stratification of $\,\overline{Jac}\,$
in terms of Gr\"obner cells for $e=0$. 
 Given a standard $M$, let
$M'\equal z^{v}M$ for $v=v_M=\de-\dim C[[z]]/M$. 
We will describe $\tilde{\la}(M')$ for the map
$M'\mapsto \tilde{I}'=\pi^0(M')=\tilde{I}(M')$, where
the latter is the lift 
of $z^{2\de}M'$ to $\C[[x,y]]$ (the inverse image).
Note that $\dim \C[[z]]/M'= \dim \tilde{C}/\tilde{I}'$.  
Generally, 
\begin{align}\label{z2dee}
\pi^e: \{M \mid z^{2\de+e}M\subset \r\}
\xrightarrow{\sim}
\{\tilde{C}^{(2\de+e)}\subset \tilde{I}
\subset \tilde{C}^{(e)} \mid P(x,y)\in \tilde{I}\}
\end{align} 
for any $\r$\~modules $M\subset \C[[z]]$ and
ideals $\tilde{I}\subset \C[[x,y]]$. One has here:
$\dim \C[[z]]/M= e+\dim \tilde{C}/\tilde{I}$
for $\pi^e(M)=\tilde{I}$. Note that standard $M$
correspond to $\tilde{I}$ 
containing some (full) lifts
of $\phi z^{2\de+e}$ for $\phi\in 1+z\C[[z]]$. 
For $\pi^{\{e=0\}}$ we have:
\begin{align}\label{z2deep}
\pi^0: \overline{Jac}=\{M'\}
\xrightarrow{\sim}
\{\tilde{C}^{(2\de)}\subset \tilde{I}'
\subset \tilde{C} \mid P(x,y)\in \tilde{I}'\}.
\end{align} 

We have $5$ standard $\Ga$-modules $\De$ (containing $0$) and 
the corresponding $5$ families of standard modules $M\subset \C[[z]]$.
Let us describe  the partitions $\tilde{\la}$
and the  diagrams $\tilde{\la}'$
for these families considered in $\overline{Jac}$.
\vfil 

{\sf Total $M_{t\!ot}=\C[[z]]$.} This  module
has $\hat{\De}\equal\De\setminus\Ga=\{1,2,5\}$. It 
was considered after Proposition \ref{Pio}. We have: $v=3$,
$\tilde{I}(M'_{t\!ot})=\{y^3,y^2x,yx^2,x^3,\ldots\}$ and 
$\tilde{\la}_{t\!ot}'=\yng(3,2,1)$. This family is just
one point.
\vskip 0.2cm
\vfil

{\sf Family 0: invertibles.} The condition  
$\hat{\De}=\emptyset$ is necessary and sufficient. These modules
are $M=\phi \r\subset \C[[z]]$, where 
$\phi=1+\al z+\be z^2+\ga z^5$ and $\al,\be,\ga\subset \C$
give the parametrization of $Jac$. Here $M'=M$ and we need to lift
$t^6 M$ to $\C[[x,y]]$. The result is $\tilde{I}'$ generated by 
the lift of $z^6\phi=z^6+\al z^7+\be z^8+\ga z^{11}$, which is 
${\vph}\equal y^2+\al yx+\be x^2+\ga x^2y$, and by 
 the lift of $(z^{4\de})$, which is 
$\tilde{C}^{(6)}=
\lan y^4, x^3, y^3x, y^2x^2, yx^3, x^4, \cdots\ran$. Note that
the latter contains $-P=x^3-y^4$, a special feature of
this example simplifying a bit the considerations.  
The cases are as follows.
\vskip 0.2cm

{\sf (0,i): $\be\neq 0,\be\neq \al^2$}.  
Then $x^2,y^2x$ belong to  $\tilde{I}_0'$,
the monomial ideal for $\tilde{I}'$. To see this compare
$y{\vph}$ and $x{\vph}$ modulo $\tilde{C}^{(6)}$.
The corresponding $\tilde{\la}'_0$ is $\yng(4,2)$.  
This sub-family is isomorphic to $\C\times \C^\ast\times \C^\ast$
as a space.

{\sf (0,ii): $\be=0, \al\neq 0$}. Then  the lowest Gr\"obner 
monomial of $y^2{\vph}= y^2+\al yx+ \ga x^2y$ is $yx$. Since 
$x^3,y^4\in 
\tilde{I}'$, we obtain:\, $\tilde{\la}'_0=\yng(4,1,1)$. As a space, 
this sub-family
is isomorphic to $\C \times \C^\ast$.

{\sf (0,iii): $\be\neq 0, \be=\al^2$}. In this case 
$y^2x\not\in \tilde{I}_0'$ as in $(i)$. However, 
$y^3\in \tilde{I}_0'$ 
due to $y{\vph}-\al x{\vph}=y^3 \mod 
 \tilde{C}^{(6)}$. Thus $\tilde{\la}'_0=\yng(3,3)\,$ and the 
corresponding space
is $\C\times \C^\ast$.

{\sf (0,iv): $\al=0=\be, \ga\neq 0$}. Here
$\tilde{\la}'_0=\yng(3,2,1)$ and the space is $\C^*$. 

{\sf (0,v): $\al=0=\be=\ga$}. We lift  $z^6 \r$,
$\tilde{\la}'_0=\yng(2,2,2)\,$, the space is one point.
\vskip 0.2cm

The families below are those containing the monomial ideals
$M_1, M_2$ and $ M_3$, considered above.
\vskip 0.2cm
\vfil

{\sf Family 1: through $M_{1}$.} Such
$M$ are  generated by
$1+\al z^2, z+\be z^2$. One has: $\hat{\De}=\{1,5\},
v=|\hat{\De}|=2$. Accordingly, the lift of  $z^8 M$ is
generated by $x^2+\al y^2x,\ y^3+\be y^2x$\, modulo
$\tilde{C}^{(6)}$. The cases are:

{\sf (1,i): $\be\neq 0$}. Then $y^2x\in \tilde{I}_0'$
and $\tilde{\la}_1'= \yng(4,2)$.  The space is 
$\C \times \C^\ast$. The Gr\"obner $f$-generators will be
$x^2-\frac{\al}{\be}y^3, y^2x+\frac{1}{\be}y^3$. 

{\sf (1,ii): $\be=0$}. Then $\tilde{\la}_1'=\yng(3,3)\,$,
 and the space
is $\C$. 
\vskip 0.1cm
\vfil

{\sf Family 2: through $M_{2}$.} Here $M$ are
generated by $1+\al z, z^2$, and  $ \hat{\De}= \{2,5\}, v=2$.
The module $z^8 M_2$ is generated by $x^2+\al y^3, y^2x$
modulo $\tilde{C}^{(6)}$. Thus we have only one subcase here:

{\sf (2,i):\ } $\tilde{\la}'=\yng(4,2)$, and 
the corresponding space is $\C$.
\vskip 0.2cm

{\sf Family 3: through $M_3$.}
The modules $M$ are generated by $1+\al z+\be z^2, z^5$, and
$\hat{\De}=\{5\}, v=1$; they  also contain $z^3,z^4,z^5\cdots$.
So the lift $\tilde{I}(M')$ of $z^7 M$ has the generators
$yx+\al x^2+\be y^3, y^2x, yx^2$ modulo $\tilde{C}^{(6)}=
\lan y^4, x^3, y^3x, y^2x^2, yx^3, x^4,\ldots \ran$. The
cases are:

{\sf (3,i): $\al\neq 0$}. Then $x^2\in \tilde{I}(M')$,
$\tilde{\la}'_3=\yng(4,2)$, and
the space is $\C\times \C^\ast$. 

{\sf (3,ii): $\al=0$}. Then $yx\in \tilde{I}(M')$,
$\tilde{\la}'_3=\yng(4,1,1)$, and
the space is $\C$. 
\vskip 0.2cm
\vfil

{\sf Summary.} The main purpose of this calculation is to
decompose $\overline{Jac}$ using the Gr\"obner cells.
E.g., the portion of $\pi^0(\overline{Jac})$ 
corresponding to $\la'=\yng(4,2)$ for $\la=\{4,2\}$
is 
$(\C\times \C^\ast\times \C^\ast)_{0,i} \cup 
(\C\times \C^\ast)_{1,i}\cup (\C)_{2,i} \cup  
(\C\times \C^\ast)_{3,i}$,
where the suffix shows the source of this contribution. 
 In the Grothendieck ring $K_0(var/\C)$, it
is $\C\times(\C-pt)^2+2 \C(\C-pt) +\C=
\C^3-2\C^2+\C+2\C^2 -2\C+\C=\C^3$. 

This is the subset in the full Gr\"obner cell $C_\la$ in
$H^{(6)}$  of the ideals $\tilde{I}$ satisfying the embeddings
\begin{align}\label{C-incl}
\tilde{C}^{(6)}=\lan x^3, y^3x,\ldots \ran\, \subset\,  
\tilde{I}\,\subset\, \tilde{C}=\lan y^2, yx, x^2,\cdots\ran.
\end{align}
Any ideals $\tilde{I}\in C_\la$ have generators:
$f_{20}=x^2+ c xy+dy^2+e y^3,\, 
f_{12}=xy^2+gy^3\in \tilde{I}$, and $f_{04}=y^4$. 
The term $y$ is missing 
here in the first generator  for {\em any\,} 
$\tilde{I}\in H^{(6)}$ 
corresponding to this diagram due to the equality $C_{01}^1=0$
from (\ref{ex-6}); see there and
Figure \ref{6-diag}. This is also granted due to
$\tilde{I}\,\subset\, \tilde{C}$.

Since $yf_{20},y^3\in \tilde{I}$, 
$x^2y+cxy^2+dy^3\in \tilde{I}$.
Using now that $x^3,xy^3\in \tilde{I}$, we obtain that
$xf_{20}\in \tilde{I}\Rightarrow cx^2y+dxy^2\in \tilde{I}$.
Thus, $(c^2-d)xy^2+cdy^3\in \tilde{I}$. Combining it with
$f_{12}$, we obtain that $cd=(c^2-d)g$ is the equation
of $\pi^0(\overline{Jac})\cap C_\la$, where $e$ is arbitrary. 

 Note that $(c\!=\!0,d\!=\!0,
f\!=\!0)$ 
is a singular point of the surface $\{cd=(c^2-d)g\}$.
 This relation reduces $\dim C_\la$ 
from $4=6-2$ to $3$. It combines in one equality
all relations above for $\al,\be$
from different subcases. Namely,
the equality $c=0=d$ (any $g,e$) corresponds to $(1,i)\&(2,i)$.
When $c^2-d\neq 0$ (any $e$), we obtain 
$(0,i)\&(3,i)$. 
\vskip 0.2cm

As another example, let us consider $\la'=\yng(3,3)\,$. 
Then the intersection of $\pi^0(\overline{Jac})$ with
the Gr\"obner cell $C_\la$  will be 
$(\C\times \C^\ast)_{0,iii}\cup (\C)_{1,ii}=\C^2 $ in the
Grothendieck ring. Now the generators of $\tilde{I}'$
subject to (\ref{C-incl}) are $x^2+cxy +dy^2+exy^2$ and $y^3$,
subject to the  relation  $d=c^2$. Also, 
$\{pt\}_{t\!ot}$ combined with 
$(\C^*)_{0,iv}$ gives $\C$ in the case of
$\,\yng(3,2,1)\,$. 

The monomiality from part (ii) of Theorem \ref{thmtop}
combined with \cite{BB} provide some a priory reasons for these 
intersections to be (topologically) affine spaces.
\vskip 0.1cm

Let us extend the first step of our
calculation to $a=3$ and $b=m$.
%\vskip 0.1cm

\begin{proposition} Let 
$\r=\C[[x\!=\!z^m,y\!=\!z^3]]$, $2\de=2(m-1)$, assuming that 
either (i)\, $m=2+3k$ or (ii)\, $m=1+3k$ for $k\ge 1$.
The partition for generic standard invertible $M$  
for the Gr\"obner decomposition of $Jac$, 
the generalized Jacobian $Jac$, under $\pi^0$ will
be denoted $\tilde{\la}_0$.

Then $\tilde{\la}_0=\{4k+2, 2k\}$ in case (i) and the
corresponding manifold is $\C^*\times \C^{\de-1}$. 
In case (ii), the partition 
is $\{4k,2k\}$. The corresponding manifold
is $\C^{\de-2}$ times  $\C^2$ minus a union of two different
$\C^1$ inside it. 
\end{proposition}

\vfil
{\it Proof.} 
We represent the generators of the invertible
modules constituting $Jac$ as 
$\phi=1+\al z+\be z^2+\ga z^4+\ldots\,.$ 
Then $2\de=2(m-1)$ is $m+3k$ for $m=2+3k$ (case $(i)$)
or $6k$ for $m=1+3k$ (case $(ii)$). Accordingly, 
$\vph=z^{2\de}\phi=xy^k+\al y^{2k+1}+\be x^2+\ga y^{2k+2}+\ldots$
in case $(i)$ 
or
$\vph=z^{2\de}\phi=y^{2k}+\al xy^{k}+\be x^2+\ga xy^{k+1}+\ldots$
in case $(ii)$. 

\vfil
{\it Case $(i)$.} The ideal  $\pi^0(z^{2\de})$ is linearly 
generated by the lift of $z^{4\de}\C[[z]]$ to $\C[[x,y]]$, which
is $\{x^2y^{2k},xy^{3k+1},x^3y^k=y^{4k+2}, x^2 y^{2k+1},\cdots\}$.
Thus it contains $y^{4k+2}$ and the first row of 
$\tilde{\la}_0'$ has at most 
$4k+2$ boxes. One has: $x\vph= x^2y^k+\al xy^{2k+1}+\be x^3+
\ga xy^{2k+2}+\cdots=
-\frac{1}{\be}(xy^k+\al y^{2k+1}+\ga y^{2k+3}+\cdots)y^k+
\al xy^{2k+1}+\be y^{m}+\ga xy^{2k+2}$. The smallest power of
$y$ here is in $xy^{2k}$, which gives that the second row is
at most with $2k$ boxes. However, $4k+2+2k=6k+2=2\de$, so
$\tilde{\la}_0$ is exactly $\{4k+2, 2k\}$ in this case. 
The corresponding manifold is $\C^*\times \C^{\de-1}$.

{\it Case $(ii)$.} The ideal  $\pi^0(z^{2\de})$,
the lift of $z^{4\de}\C[[z]]$, is now 
linearly generated by
$\{y^{4k}, xy^{3k}, x^2y^{2k}, xy^{3k+1},\cdots\}$.
So it contains $y^{4k}$ and the first row of 
$\tilde{\la}_0'$ has at most 
$4k$ boxes. One has: $x\vph= 
xy^{2k}+\al x^2y^{k}+\be x^3+\ga x^2y^{k+1}+\ldots=
xy^{2k}
-\frac{\al}{\be}(y^{2k}+\al xy^{k}+\ga xy^{k+1}+\ldots+\cdots)y^k+
+\be y^{m}+\be y^{m}+\ga x^2y^{k+1}$. The term with
the smallest power of
$y$ is now $(1-\frac{\al^2}{\be})xy^{2k}$, 
which gives that the second row is
at most with $2k$ boxes provided that $\be\neq \al^2$.
 However, $4k+2k=6k=2\de$, so
$\tilde{\la}_0$ is exactly $\{4k, 2k\}$ in this case. 
The corresponding manifold is isomorphic
to $\C^2\setminus (\C\times pt\,\cup\, 
pt\times \C)$ multiplied by $\C^{\de-2}$. \sq

\vfil
Without going into all detail,
let us provide the Young diagram $\tilde{\la}'$ for generic
points of $Jac$ for the ring $\r=\C[[x\!=\!z^5,y\!=\!z^4]]$. 
It is with $2\de= 12$ boxes:\, $\,\yng(6,4,2)\,$. 
First of all,
$y^6=t^{24}=t^{4\de}$ and $y^6\in \tilde{I}$. Setting
$\vph=z^{2\de}\phi=y^3+\al xy^2+\be x^2 y+\ga x^3 +\de x^2y^2+
\ep x^3y+\nu x^3y^2$, one assumes that $\ga\neq 0$,
which gives that $x^3$
belongs to the corresponding monomial $\tilde{I}_0$.
Considering  $x\vph$, we obtain that 
$(\al-\frac{\be^2}{\ga})x^2y^2$
belongs to $\tilde{I}_0$, so we assume next that $\be^2 \neq \al\ga$
and obtain that $x^2y^2\in \tilde{I}_0$.
Similarly, $x^2\vph$, $xy\vph$, and $y^2\vph$
 give that  
$x^2y^3+\al x^3 y^2$, $xy^4+\al x^2y^3+\be x^3y^2$,
$xy^4+(\be-\al^2) x^3y^2$,
 and
$\al x y^4+(\ga-\be\al) x^3y^2+y^5$ belong to 
$\tilde{I}$. Combining the latter two elements, we obtain
that $xy^4\in \tilde{I}$ if $\ga\neq 2\al\be +\al^3$.
\vfil

For arbitrary relatively prime $b>a>1$ and 
$\r=\C[[x\!=\!z^b,y\!=\!z^a]]$, the Young diagrams
for $\pi^0(\overline{Jac})$ 
are of order $2\de=(a-1)(b-1)$ and belong to the
rectangle $a\times \kappa$ for
$\kappa=Ceiling[4\de/a]$. We use that the
corresponding ideals contain $\pi^0(z^{2\de}\C[[z]])$,
so $y^\kappa$ belongs to all of them.
The polynomial $x^a-y^b$ gives that the second dimension
is $a$. 

Recall, that $\tilde{I}\subset 
\pi^0(\C[[z]])$, i.e. $\tilde{\la}'$
contains the diagram $\la'_{t\!ot}$ 
for the lift of $z^{2\de}\C[[z]]$ to
$\C[[x,y]]$. For instance, 
the minimal possible number of rows in
$\tilde{\la}'$ is $Ceiling[2\de/b]=a-1$.

In examples, all such diagrams can be
obtained from some modules $M$ in $\overline{Jac}$.
Their total number of diagrams with $2\de$ boxes in
the rectangle $a\times \kappa$ equals the coefficient of 
$q^{2\de}$ in $q^{\binom{a}{2}}\binom{a+\kappa}{a}_q$
for the $q$-binomial coefficients $\binom{n}{m}_q$.
We need to diminish it by the number of such diagrams with
$a-2$ rows.
Generally, this difference is greater than the rational
slope Catalan number $\frac{1}{a+b}\binom{a+b}{a}$, which
gives the number of Piontkowski cells. Recall that each 
such cell generally results in several diagrams. 
\vfil

If all such partitions
$\tilde{\la}$ 
can be obtained this way, then 
the connectivity of $\cup_{\tilde{\la}}\, C_{\tilde{\la}}\,\subset 
H^{(2\de)}$ follows from that for $\overline{Jac}$, which can be
of independent interest.  
Actually only the diagrams with the minimal possible number
of rows, which is $a-1$, are sufficient here to check, and 
we can use different
orderings of $\{x^m y^n\}$, not only the one with 
$\{x^\infty\!<\!y\}$. 

\vskip 0.2cm
\vfil

{\sf Some perspectives.}
Summarizing, we presented  $\overline{Jac}$ above
as the union of
the intersections of its $\pi^0$\~image with the proper
Gr\"obner cells. The resulting
intersections are homeomorphic to $\C^3, \C^2, \C^2, \C, pt$, 
i.e. the same 
as for the Piontkowski decomposition with respect to
$\De(M)$. The list of cells must be the same because their 
multiplicities are Betti numbers of $\overline{Jac}$.
Another, similar, approach is to consider the filtration 
of $\overline{Jac}$ in terms of closures of the sets
of invertibles inside monomial standard $M$ and the corresponding
strata. Combined with the Lusztig-Smelt-Piontkowski cells, it
gives a justification of the super-duality of the
{\em motivic superpolynomials\,} for quasi-homogeneous $\r$, to
be discussed elsewhere.  
The Gr\"obner decomposition of $\pi^0(Jac)$ (Family 0: invertibles) is 
interesting in its own right, with possible relations to \cite{MY,MS}.
\vfil

Generally,  we can decompose
the images of $\pi^0(\overline{Jac})$ 
due to the presentation from (\ref{z2deep})
using \cite{BB} and the methods based on the 
{\em stable envelopes\,}  from  \cite{MO}. This is
for the action of $\C^*$ for quasi-homogeneous singularities.
Actually, it suffices to assume here
that $\tilde{C}_\r$ and the image of $P(x,y)$ modulo
$\tilde{C}^{(2\de)}_\r$ are invariant with respect to
the action $\,x\mapsto u^b x,\, y\mapsto u^ay$ for $u\in \C^\ast$.
The method from \cite{BB} provides then that the
corresponding "cells" will be affine spaces. 
Actually, only the $C^\ast$-invariance
of $\overline{Jac}$ and $ Jac^\bullet$ is sufficient. For instance,
Theorem \ref{thmtop} (its end) provides this invariance when
$\Ga$ has $2$ generators or $3$ with $p=1$. 
\vskip 0.2cm 
\vfil

{\sf Using different orderings.}
There is another approach, adjusted directly
to quasi-homogeneous plane curve singularities. Its aim is
to eliminate the non-trivial combinatorics of the 
Gr\"obner decomposition of $Jac$ (Family 0) and other 
Lusztig-Smelt-Piontkowski cells.
Given $a,b\in \Z_+$ such that $gcd(a,b)=1$,
one introduces the {\it valuation\,} $\nu(y)=a$, $\nu(x)=b$
in $\C[[x,y]]$. Then the {\it weighted Gr\"obner basis\,} and
the corresponding $\tilde{I}_0$ are defined when we first
order monomials with respect to $\nu$ (the smallest $\nu$ first),
and then with respect to our usual Gr\"obner ordering  
$\{x^\infty\!<\!y\}$ if their $\nu$ coincide.
This leads to a variant of the {\em wall crossing\,}, where
the ratio $a/b$ serves as the {\em stability condition\,}. 
We will not discuss here the corresponding theory; the 
approach from \cite{BB}-\cite{MO} can be used.  Let us give
one example.
\vskip 0.2cm

We will calculate all  $\tilde{\la}_\nu'$ for  
$\r=\C[[y=z^a,x=z^b]]$, 
and $Jac$ (Family 0) for the example above.
Let $a=3, b=4$, $\nu=\nu_z$.
The smallest monomial in 
${\vph}=y^2+\al yx+\be x^2+\ga x^2y$ will be then always $y^2$, 
 and the corresponding
{\em weighted\,} $\tilde{\la}_\nu'$ becomes $\yng(2,2,2)\,$ for 
{\em any\,}
$\al,\be,\ga$. Only {\em monomial\,} $M$ are sufficient
to consider. Let us list all $\tilde{\la}'_\nu$ for the
corresponding families:

$$ (t\!ot): \yng(3,2,1),\, (0): \yng(2,2,2)\,,\,
(1): \yng(3,3)\,,\, (2): \yng(4,2),\,
(3): \yng(4,1,1).
$$ 

This is relatively straightforward for {\em arbitrary\,} 
quasi-homogeneous singularities. We represent $2\de=(a-1)(b-1)=
au+bv$ for $u,v\in \Z_+$. In the case of  "the big cell" $Jac$, 
the weighted $\tilde{\la}_\nu'$ will be the rectangle 
$a\times (b+u)$ without the corner 
$\{\{i,j\} \mid i\ge v\,\&\, j\ge u\}$. 
The corresponding monomial ideal is generated by  $x^a, y^{b+u},
x^v y^u$.
Note that the number of boxes of this diagram is\,
$a(b+u)-(a-v)b=au+bv=2\de$, as it is supposed to be. The number
of lines is always no greater than $a$, since $\pm\, P(x,y)=x^a+\ldots$ 
belongs to all ideals for $\pi^0(\overline{Jac})$. The smallest number
of columns is $(b-1)-Floor[\frac{b-1}{a}]$, which
is for $M_0$.
We see that every Piontkowski cell naturally maps
to the corresponding (single) Gr\"obner cell under the
ordering based on such $\nu$. 

Thus the Piontkowski decomposition can be generally
seen as a particular case of that based on \cite{BB}, the
theory of {\em stable envelopes\,}, and "localization theory". This
is for quasi-homogeneous plane curve singularities, but can be
possibly extended to those
in Theorem \ref{thmtop}. In Proposition \ref{MIast}, this
gives the description of the scheme $\tilde{Jac}[\De]$ entirely 
in terms
of $\Ga$ (for such singularities).

%\medskip
\bibliographystyle{unsrt}

%\vfill\eject

%%%%%!!!!! APPENDICES-COMMENT

\end{document}